\documentclass[a4,11pt]{article}
\usepackage{amsmath}
\usepackage{graphics}
\usepackage{latexsym}
\usepackage{amsfonts}
\numberwithin{equation}{section}

\newtheorem{corollary}{Corollary}
\newtheorem{lemma}{Lemma}[section]
\newtheorem{remark}{Remark}

\newtheorem{assumption}{Assumption}
\def\qed{ \ \vrule width.2cm height.2cm depth0cm\smallskip}
\def \et{\eta}
\def \ind{1\!\!1}
\def \x{X^{t,x}}
\def \udl{\underline}
\def\fin{\infty}
\def \ij {i\in \cJ}
\def \fl {\forall}
\def \orw{\overrightarrow}
\newcommand{\integ}[2]{\displaystyle \int_{#1}^{#2}}

\newcommand{\eps}{\varepsilon}

\newcommand{\brm}{\begin{rem}}
\newcommand{\ermq}{\end{rem}}

\newcommand{\ba}{\begin{array}}
\newcommand{\ea}{\end{array}}
\newcommand{\be}{\begin{equation}}
\newcommand{\ee}{\end{equation}}
\newcommand{\bea}{\begin{eqnarray}}
\newcommand{\eea}{\end{eqnarray}}
\newcommand{\beaa}{\begin{eqnarray*}}
\newcommand{\eeaa}{\end{eqnarray*}}

\def \R{I\!\!R}
\def \E{\mathbb{E}}

\def\a{\alpha}
\def\b{\beta}
\def\g{\gamma}
\def\d{\delta}
\def\e{\varepsilon}

\def\l{\lambda}

\def\n{\nu}
\def\si{\sigma}
\def\t{\tau}

\def\vro{\varrho}
\def\D{\Delta}

\def\cB{{\cal B}}
\def\cC{{\cal C}}
\def\cD{{\cal D}}

\def\cF{{\cal F}}

\def\cH{{\cal H}}

\def\cJ{{\cal J}}

\def\cL{{\cal L}}

\def\cS{{\cal S}}

\def\no{\noindent}

\def\ms{\medskip}
\def\bs{\bigskip}
\def\q{\quad}
\def\qq{\qquad}

\def\bF{{\bf F}}

\def\qed{ \hfill \vrule width.25cm height.25cm depth0cm\smallskip}

\newcommand{\basa}{\begin{assumption}}
\newcommand{\easa}{\end{assumption}}

\newcommand{\bas}{\begin{assum}}
\newcommand{\eas}{\end{assum}}

\def\limsup{\mathop{\overline{\rm lim}}}
\def\liminf{\mathop{\underline{\rm lim}}}

\def\esssup{\mathop{\rm esssup}}
\def\argmax{\mathop{\rm argmax}}

\def\vp {\varphi}
\def\dis{\displaystyle}

\def\bF{{\bf F}}

\def \P{\mathbb{P}}
\newtheorem{thm}{Theorem}[section]

\newtheorem{prop}[thm]{Proposition}
\newtheorem{rem}[thm]{Remark}

\newtheorem{assum}[thm]{Assumption}

\def\mb{\mbox}
\title{Viscosity Solutions of Systems of PDEs with Interconnected Obstacles and Multi-Modes Switching Problem}
\author{S.Hamad\`ene\thanks{Universit\'e du Maine, LMM, Avenue Olivier Messiaen, 72085 Le Mans, Cedex 9,
France ; hamadene@univ-lemans.fr}\,\,\, and \, M.A.
Morlais\thanks{Same address ; e-mail:
Marie$_-$Amelie.Morlais@univ-lemans.fr} }
\begin{document}
\date{}
\maketitle
\newtheorem{theo}{Theorem}
\newtheorem{problem}{Problem}
\newtheorem{pro}{Proposition}
\newtheorem{cor}{Corollary}
\newtheorem{axiom}{Definition}
\newcommand{\rw}{\rightarrow}
\def \R{\mathbb{R}}
\def \cadlag {{c\`adl\`ag}~}
\def \esssup {\mbox{ess sup}}
\def \esp {[0,T]\times \R^k}

\begin{abstract}

This paper deals with existence and uniqueness, in viscosity sense,
of a solution for a system of $m$ variational partial differential
inequalities with inter-connected obstacles. A particular case of
this system is the deterministic version of the Verification Theorem
of the Markovian optimal $m$-states switching problem. The switching
cost functions are arbitrary. This problem is connected with the
valuation of a power plant in the energy market. The main tool is
the notion of systems of reflected backward stochastic differential
equations with oblique reflection.
\end{abstract}

\no{\bf AMS Classification subjects}: 60G40 ; 62P20 ; 91B99 ; 91B28
; 35B37 ; 49L25.
\medskip

\no {$\bf Keywords$}: Real options; Backward stochastic differential
equations; Snell envelope; Stopping times ; Switching; Viscosity
solution of PDEs; Variational inequalities.

\section {Introduction}

The multi-modes switching problem is by now well documented both in
the economics or mathematical literatures (see e.g.
\cite{Brennanschwarz, carmonaludkovski_1, carmonaludkovski_2, Dixit,
Djehicheetal09,Hamadene_elasri, Hamjeanblanc, hamzhang, hutang,
Porchet-touzi, tangyong, dz}, etc. and the references therein). The
pioneering work of Brennan and Schwarz \cite{Brennanschwarz} deals
with a two-modes switching problem describing the life cycle of an
investment in the natural resource industry. A major switching
problem of interest is related to the energy market. Actually let us
consider a power plant which has several modes of production and
which is put in a specific mode according to its profitability which
depends on the electricity price in the market. The manager of the
plant aims at maximizing its global profit. For this objective, she
implements an optimal strategy $\delta^{*}$ which is a pair of two
sequences $(\tau_{k})_{k\geq 1}$ and $(\xi_{k})_{k\geq 1}$
describing respectively the optimal successive switching times and
modes. However switching the plant from one mode to another is not
free and generates expenditures and, on the other hand, when the
plant is in a specific mode it provides a profit which depends on
that mode. \ms

It is well-known that the optimal switching problem is related to
systems of reflected backward stochastic differential equations
(BSDEs for short) with inter-connected obstacles or oblique
reflection (see e.g. \cite{carmonaludkovski_1, Djehicheetal09,
Hamjeanblanc, hamzhang, hutang, Porchet-touzi}) of the following
type: $\forall i\in \cJ:=\{1,...,m\}$, \be \label{multidimrbsde0}
\left\{
\begin{array}{l}
Y^i_s=h^i_T(\omega)+\int_s^Tf_i(\omega,r)dr+K^i_T-K^i_s-\int_s^TZ^i_rdB_r,\,\,\forall\,\,s\leq
T\\
Y^i_s\geq \max_{j\in {\cal J}^{-i}}\{Y^j_s-g_{ij}(\omega,s)\},
\,\,\forall s\leq T\\
\int_0^T(Y^i_s-\max_{j\in {\cal
J}^{-i}}\{Y^j_s-g_{ij}(\omega,s))dK^i_s=0,
\end{array}\right.
\ee where ${\cal J}^{-i}={\cal J}-\{i\}$. Actually it is shown in
the aforementioned papers that $Y^i_0$ is the optimal profit if the
plant is in mode $i$ at $t=0$, i.e., $\forall \ij$,
\begin{equation}\label{eq:optimalvaluefunction}
\dis{ \quad Y_{0}^{i} = \sup_{\delta = (\tau_{k}, \xi_{k})_{k\geq 0}
\in \mathcal{A}_{0}^i} \E[h^{\d}(T, \omega)+ \int_{t}^{T}f_{\d}(s,
\omega)ds -\sum_{k \geq 0} g_{\xi_{k},
\xi_{k+1}}(\tau_{k+1})\mathbf{1}_{[\tau_{k+1}<T]} ]}
\end{equation}
where, $\mathcal{A}_{0}^i$ stands for the set of all admissible
strategies starting from mode $i$ at time $t=0$, $f_{\d}$ the
instantaneous profit per unit of time when $\d$ is implemented,
$g_{\xi_{k}, \xi_{k+1}}$ is the switching cost from mode $\xi_k$ to
mode $\xi_{k+1}$ ($\xi_0=i,\,\,\tau_0=0$) and finally $h^\d(T)$ is
the terminal profit under $\d$.

Additionally, if we are given the solution $(Y^{i})_{i =1,\cdots, m}$ then the optimal strategy $\delta^{*} = (\xi_{k}^*, \tau_{k}^*)$ is uniquely characterized as follows: if we set $ \tau_{1}^{*} =0,$ and $ \xi_{1}^{*} = i $ then for any $k \ge 1$
 \be   \left\{ \begin{array}{l}
 \tau_{k+1}^* =\inf \{t >\tau_{k}^*,\; Y_{t}^{\xi_{k}^*} = \max_{ j\in {\cal J}^{-\xi_{k}^*}}( Y^j_{t} - g_{\xi_{k}^*, j}(\omega,t))  \}\wedge T,\\
 \\
 \xi_{k+1}^{*} \in \argmax \{j, \;  Y_{t}^{\xi_{k}^*} = \max_{ j\in {\cal J}^{-\xi_{k}^*}} (Y^j_{t} - g_{\xi_{k}^*, j}(\omega,t) ) \}.
 \end{array}
  \right.   \ee
 In the case of Markovian setting and when randomness stems from an
exogenous  It\^o diffusion process $(X^{t,x}_s)_{s\leq T}$,
$(t,x)\in \esp$, of the following type:
$$
dX^{t,x}_s=b(s,X^{t,x}_s)ds+\sigma(s,X^{t,x}_s)dB_s, \,\,s\in
[t,T]\,;X^{t,x}_s=x \mbox{ for }s\leq t, $$that is to say, the
processes $f_i(s,\omega)$, $h_i(s,\omega)$ and $g_{ij}(s,\omega)$
are deterministic functions of $(s,X_s^{t,x}(\omega))$, the optimal
switching problem is also related to the following system of
variational inequalities with inter-connected obstacles: $\forall
\,\,i\in {\cal J}$ \be \label{sysvi0}  \left\{
\begin{array}{l}
\min\left \{v_i(t,x)- \max\limits_{j\in{\cal
J}^{-i}}(-g_{ij}(t,x)+v_j(t,x)), -\partial_tv_i(t,x)- {\cal
L}v_i(t,x)-f_i(t,x)\right\}=0;\\
v_i(T,x)=h_i(x)
\end{array}\right.
\ee where $\cL$ is the infinitesimal generator associated with
$X^{t,x}$. The process $X^{t,x}$ can be the electricity price in the
market or the dynamics of factors which determine that price. System
(\ref{sysvi0}) is the deterministic version of the verification
theorem of the optimal switching problem.

Actually in \cite{Hamadene_elasri}, the authors have proved that if,
mainly, the switching costs satisfy $g_{ij}(t,x)\geq \g_0>0$ then
the system (\ref{sysvi0}) has a unique continuous solution
$(v^1,\dots,v^m)$ in viscosity sense and the following relationship
holds true:
$$
\forall \ij, (t,x)\in \esp, \forall s\in [t,T],
\,\,Y^i_s=v^i(s,X^{t,x}_s).$$ where $(Y^i)_{\ij}$ are the processes
solution of (\ref{multidimrbsde0}) with $f_i$, $g_{ij}$ and $h_i$ which
are deterministic functions of $(s,X_s^{t,x})$.

So the main objective of this paper is to deal with system
(\ref{sysvi0}) in its general setting, i.e., to study the existence
and uniqueness of a solution in viscosity sense for the following:
$\forall \,\,i\in {\cal J}$ \be \label{sysvi01} \left\{
\begin{array}{l}
\min\left \{v_i(t,x)- \max\limits_{j\in{\cal
J}^{-i}}(-g_{ij}(t,x)+v_j(t,x)),\right.\\\left.
\qquad\qquad-\partial_tv_i(t,x)- {\cal
L}v_i(t,x)-f_i(t,x, v^1(t,x), \dots,v^m(t,x), \sigma^\top(t,x)D_xv^i(t,x))\right\}=0\,\,;\\
v_i(T,x)=h_i(x).
\end{array}\right.
\ee As previously pointed out, system (\ref{sysvi01}) is connected
with the following system of reflected BSDEs with oblique reflection
: \be \label{multidimrbsde1} \left\{
\begin{array}{l}
Y^{i;t,x}_s=h^i(X^{t,x}_T)+\int_s^Tf_i(r,X^{t,x}_r,
Y^{1;t,x}_r,\dots,Y^{m;t,x}_r,Z^{i;t,x}_r)dr+K^{i;t,x}_T-K^{i;t,x}_s-\int_s^TZ^{i;t,x}_rdB_r,\,\,s\leq
T;\\
Y^{i;t,x}_s\geq \max_{j\in {\cal
J}^{-i}}\{Y^{j;t,x}_s-g_{ij}(s,X^{t,x}_s)\},
\, s\leq T;\\
\int_0^T(Y^{i;t,x}_s-\max_{j\in {\cal
J}^{-i}}\{Y^{j;t,x}_s-g_{ij}(s,X^{t,x}_s))dK^{i;t,x}_s=0.
\end{array}\right.
\ee This system has already been studied in
\cite{chassagneuxetal2010,hamzhang} in this general form. \ms

There are at least two motivations for considering this problem. The
first one is to extend, as much as possible, the Feynman-Kac's
representation of the solution of (\ref{multidimrbsde1}) via the
solution of (\ref{sysvi01}) and the process $X^{t,x}$. This issue is
very important if we are willing to consider the numerical study of (\ref{multidimrbsde1}) especially in connection with the
pricing of gas options in the energy market (one can see e.g.
\cite{marieb} for more details). The second one is that solving this problem is a step toward the study of zero-sum
switching games which are encountered especially in the carbon
market by energy firms (see e.g. \cite{marieb} for more details).
\ms

The novelty of this paper lies in the fact that we investigate both
existence and uniqueness of a continuous viscosity
solution of (\ref{sysvi01}) under the following relaxed hypotheses:\\
\noindent (i) the switching costs are non-negative and satisfy the no
free loop condition (see [H3]-(ii) below)

\noindent (ii) for any $\ij$, either :

(a) $\forall\,\,k\in\cJ^{-i}$,   the mappings $y_k\rw
f_i(t,x,y^1,\dots,y^{k-1},y^k,y^{k+1},\dots,y^m,z)$ are
non-decreasing ;

\noindent or

(b) $\forall\,\,k\in\cJ^{-i}$,   the mappings $y_k\rw
f_i(t,x,y^1,\dots,y^{k-1},y^k,y^{k+1},\dots,y^m,z)$ are
non-increasing. \ms

In both cases we show existence of a solution for (\ref{sysvi01})
while we have been able to show uniqueness only in the case
(ii)-(a). \ms

The closest paper to ours is the one by Elie-Kharroubi
\cite{eliekharroubi} where the authors deal also with the
representation of solution of (\ref{multidimrbsde1}) by the
viscosity solution of (\ref{sysvi01}). However their approach, based
on the minimal solutions of constrained BSDEs with jumps introduced
in \cite{MaKharroubi_2010}, is not very satisfactory since it
induces assumptions on the data of the problem which are either not
natural or difficult to verify in practice. Finally note that there
are also works related to viscosity solutions of the switching
problem but their settings and/or approaches are not the same as
ours \cite{djehicheetal2,brunoboucha, hutang,tangyong}, etc. \ms

This paper is organized as follows:

In Section 2 we collect the main assumptions on the data of the
problem and we define the notion of a viscosity solution for the
system (\ref{sysvi01}). In Section 3, we once more introduce the
switching problem and provide some results related to solutions of
systems of reflected BSDEs with oblique reflection which are rather
new since they are obtained under weaker conditions than the ones of
the literature on the subject (see e.g. \cite{chassagneuxetal2010,
hamzhang}, etc.). In particular, we mainly deal with the non free
loop property [H3]-(ii) on the switching costs. Those results are
basic to deal with the main purpose of this work. In Section 4, we
provide a comparison result between sub-solutions and
super-solutions of the system (\ref{sysvi01}) in the case when for
any $i\in \cJ$ the function $f_i$ depends on $(y^i)_{\ij}$ only
through $y^i$. We then show, in this specific framework of functions
$f_i$, that system (\ref{sysvi01}) has a unique continuous solution
$(v^i(t,x))_{\ij}$ which is moreover of polynomial growth. As a
by-product, we provide a probabilistic representation for the
solution of (\ref{multidimrbsde1}) via the deterministic continuous
functions $(v^i(t,x))_{\ij}$ and the process $X^{t,x}$. In Section
5, we deal with the general framework under mainly conditions
(i)-(ii) above. Using the results of Sections 3 and 4, we construct
in each case an approximating scheme which is convergent and whose
limit is a solution in viscosity sense for system (\ref{sysvi01}).
Finally under condition (ii)-(a), we show that system
(\ref{sysvi01}) satisfies the comparison property between sub- and
super-solutions. Thus under conditions (ii)-(a) the solution of
system (\ref{sysvi01}) is unique. \qed
\section{Assumptions and problem formulation}
Let $T$ (resp. $k$) be a fixed real (resp. integer) positive
constant, let $\cJ:=\{1,\dots,m\}$ and let us consider the following
functions: for $i,j\in \cJ$,
$$
\begin{array}{l}
b:(t,x)\in [0,T]\times \R^k\mapsto b(t,x)\in \R^{k};\\
\sigma:(t,x)\in [0,T]\times \R^k\mapsto \sigma(t,x)\in \R^{k\times
d};\\ f_i: (t,x,y^1,...,y^m,z)\in [0,T]\times \R^{k+m+d}\mapsto
f_i(t,x,y^1,...,y^m,z)\,\,\in \R\,;\\
g_{ij}:(t,x)\in [0,T]\times\R^k\mapsto g_{ij}(t,x)\in \R\,(i\neq j);\\
h_i: x\in \R^k\mapsto  h_i(x)\in \R.\\
\end{array}$$

Next let $\phi: (t,x)\in \esp\mapsto \phi(t,x)\in \R$ be a function.
It is called of {\it polynomial growth} if there exist two non
negative real constants $C$ and $\g$ such that:$$|\phi(t,x)|\leq
C(1+|x|^\g),\,\,\forall (t,x)\in \esp.$$  Throughout this paper, we denote by $\Pi^g$ the class of functions with
polynomial growth and by $\cC^{1,2}(\esp)$ (or simply
$\cC^{1,2}$) the set of functions defined on $\esp$ with values in
$\R$ which are $C^1$ in $t$ and $C^2$ in $x$. $\Box$ \ms

We now consider the following assumptions: \ms

\no {\bf [H1]}: The functions $b$ and $\sigma$ are jointly
continuous and of linear growth in $(t,x)$, and Lipschitz continuous
w.r.t. $x$, i.e., there exists a constant $C\geq 0$ such that for
any $t\in [0, T]$ and $x, x'\in \R^k$\be \label{regbs1}|b(t,x)|+
|\sigma(t,x)|\leq C(1+|x|) \quad \mbox{ and } \quad
|\sigma(t,x)-\sigma(t,x')|+|b(t,x)-b(t,x')|\leq C|x-x'|.\ee

Throughout this paper we assume that assumption [H1] holds. \\

\ms
\no {\bf [H2]}: for $i\in {\cal J}$, $f_i$ satisfies:

(i) $(t,x)\mapsto f_i(t,x,y^1,...,y^m,z)$ is continuous uniformly
w.r.t. $(\overrightarrow{y},z):=(y^1,...,y^m,z)$;

(ii) $f_i$ is uniformly Lipschitz continuous with respect to
$(\overrightarrow{y},z):=(y^1,...,y^m,z)$, i.e., for some $C\geq 0$,
$$|f_i(t,x,y^1,...,y^m,z)-f_i(t,x,\bar y^1,...,\bar y^m,\bar z)|\leq C(|y^1
-\bar y^1|+\dots+|y^m-\bar y^m|+|z-\bar z|)\,\,;$$

(iii) the mapping $(t,x)\mapsto f_i(t,x,0,\dots,0)$ is ${\cal
B}([0,T]\times \R^k)$-measurable and of polynomial growth i.e. it
belongs to $\Pi^g$;\ms

(iv) \underline{\it{Monotonicity}}: $\forall i\in {\cal J}$, for any
$k\in {\cal J}^{-i}$, the mapping $y_k\in \R\mapsto
f_i(t,x,y_1,...,y_{k-1},y_k,y_{k+1},...,y_m)$  is non-decreasing
whenever the other components
$(t,x,y_1,...,y_{k-1},y_{k+1},...,y_m)$ are fixed; \bs

\noindent {\bf [H3]}: (i) $g_{ij}$ is jointly continuous in $(t,x)$,
non-negative, i.e., $g_{ij}(t,x)\geq 0$, $\forall (t,x)\in
[0,T]\times \R^k$ and belongs to $\Pi^g$;

(ii) {\it The non-free loop property}: for any $(t,x)\in [0,T]\times
\R^k$ and for any sequence of indices $i_1,...,i_k $ such that
$i_1=i_k$ and $card\{i_1,...,i_k\}=k-1$ we have:
$$g_{i_1 i_2}(t,x)+g_{i_2 i_3}(t,x)+\dots+g_{i_{k-1} i_k}(t,x)+g_{i_k i_1}(t,x)>0,
\,\,\forall (t,x)\in [0,T]\times \R^k.$$ As a convention we assume
hereafter that $g_{ii}(t,x)=0$ for any $(t,x)\in \esp$ and $\ij$.\ms

\no {\bf [H4]}: $h_i$ is continuous, belongs to $\Pi^g$ and
satisfies:
$$\forall x\in\R,\,\,
h_i(x)\geq \max_{j\in {\cal J}^{-i}}(h_j(x)
-g_{ij}((T,x)).\,\,\,\,\qed$$

Next let us introduce the following infinitesimal generator \be \cL
\varphi(t,x)
=\frac{1}{2}Tr[(\sigma.\sigma^\top)(t,x)D_{xx}^2\varphi(t,x)]+b(t,x)^\top.D_x\varphi(t,x)
\ee for a function $\varphi$ which belongs to $\cC^{1,2}([0,T]\times
\R^k;\R)$ ($Tr$ is the trace of a symmetric matrix and $(.)^\top$
stands for the transpose). It is associated with a stochastic
process which we will describe precisely below. \qed
\medskip

In this paper we are concerned with the existence and uniqueness in
viscosity sense of the solution $(v^1,\dots,v^m): (t,x)\in
[0,T]\times \R^k\mapsto (v^1(t,x),\dots,v^m(t,x))\in \R^{m} $ of the
following system of $m$ partial differential equations with
inter-connected obstacles:  $\forall \,\,i\in {\cal J}$ \be
\label{sysvi1}
  \left\{
\begin{array}{l}
\min\left \{v_i(t,x)- \max\limits_{j\in{\cal
J}^{-i}}(-g_{ij}(t,x)+v_j(t,x))\right.;\\ \qquad \qquad
\left.-\partial_tv_i(t,x)- {\cal
L}v_i(t,x)-f_i(t,x,v^1(t,x),...,v^m(t,x),\sigma^\top(t,x).D_xv^i(t,x))\right\}=0;\\
v_i(T,x)=h_i(x).\qed
\end{array}\right.
\ee

To proceed we will precise the notion of a viscosity solution of the
system (\ref{sysvi1}). It will be done in terms of sub- and
super-jets. So for any locally bounded function $u:\,(t,x)\in
[0,T]\times \R^k\mapsto u(t,x)\in \R$, we define its lower
semicontinuous ($lsc$ for short) envelope $u_*$, and upper
semicontinuous ($usc$ for short) envelope $u^*$ in the following
way: $$u_*(t,x)=\liminf_{(t',x')\rightarrow (t,x),\,t'<T}u(t',x')
\mbox{ and }u^*(t,x)=\limsup_{(t',x')\rightarrow
(t,x),\,t'<T}u(t',x').$$

\begin{axiom}: Subjects and superjets\\\ms
$(i)$ For a function $u:[0,T]\times \R^k\rightarrow \R$, lsc (resp.
usc), we denote by $J^-u(t,x)$ the parabolic subjet (resp.
$J^+u(t,x)$ the parabolic superjet) of $u$ at $(t,x)\in [0,T]\times
\R^k$, as the set of triples $(p,q,M)\in \R\times \R^k\times
\mathbb{S}^k$ satisfying \beaa u(t',x')& \geq \,(\mbox{ resp. }
\leq)\,\,u(t,x)+p(t'-t)+<q,x'-x>+\\&\qquad
\frac{1}{2}<x'-x,M(x'-x)>+o(|t'-t|+|x'-x|^2)\eeaa where
$\mathbb{S}^k$ is the set of symmetric real matrices of dimension
$k$. \bs

\no$(ii)$ For a function $u:[0,T]\times \R^k\rightarrow \R$, lsc
(resp. usc), we denote by $\bar J^-u(t,x)$ the parabolic limiting
subjet (resp. $\bar J^+u(t,x)$ the parabolic limiting superjet) of
$u$ at $(t,x)\in [0,T]\times \R^k$, as the set of triples
$(p,q,M)\in \R\times \R^k\times \mathbb{S}^k$ such that: \beaa
(p,q,M)=\lim_n
(p_n,q_n,M_n), & (t,x)=\lim_n (t_n,x_n)\\
\mbox{ with }(p_n,q_n,M_n)\in J^-u(t_n,x_n) \,\,(\mbox{resp.
}J^+u(t_n,x_n))& \mbox{and }u(t,x)=\lim_n u(t_n,x_n). \eeaa
\end{axiom}

We now give the definition of a viscosity solution for the system of
PDE equations with oblique reflection (\ref{sysvi1}).

\begin{axiom}: Viscosity solution to (\ref{sysvi1})\\
(i) A function $(v_1,...,v_m):[0,T]\times \R^k\rightarrow \R^m$ such
that for any $i\in \cJ$, $v_i$ is lsc (resp. usc), is called a
viscosity supersolution (resp. subsolution) to (\ref{sysvi1}) if for
any $i\in \cJ$, for any $(t,x)\in [0,T)\times \R^k$ and any
$(p,q,M)\in \bar J^-v_i(t,x)$ (resp. $\bar J^+v^i(t,x)$) we have:
$$\left\{\begin{array}{l}\min\left\{v_i(t,x)- \max\limits_{j\in{\cal
J}^{-i}}(-g_{ij}(t,x)+v_j(t,x));\right.\\\left.\qquad
-p-b(t,x)^\top.q-\frac{1}{2}
Tr[(\sigma\sigma^\top)(t,x)M]-f_i(t,x,v_1,\dots,v_m,\sigma
(t,x)^\top.M)\right\}\geq 0 \,\,(\mbox{resp. }\leq 0);\\v_i(T,x)\geq
\, (\mbox{resp.} \leq)\, h_i(x).\end{array}\right.$$ (ii) A locally
bounded function $(v_1,...,v_m):[0,T]\times \R^k\rightarrow \R^m$ is
called a viscosity solution to (\ref{sysvi1}) if
$({v_1}_*,...,{v_m}_*)$ (resp. $(v_1^*,...,v_m^*)$) is a viscosity
supersolution (resp. subsolution) of (\ref{sysvi1}).\qed
\end{axiom}

As pointed out previously we will show that system (\ref{sysvi1})
has a unique solution in viscosity sense. A particular case of this
system is the deterministic version of the optimal $m$-states
switching problem which is well documented e.g. in
\cite{Djehicheetal09, hutang} and which we will describe in the next
section.
\section{The  optimal $m$-states  switching problem}
\subsection{Setting of the problem}
Let $(\Omega, {\cal F}, \P)$ be a fixed probability space on which
is defined a standard $d$-dimensional Brownian motion
$B=(B_t)_{0\leq t\leq T}$ whose natural filtration is
$(\cF_t^0:=\sigma \{B_s, s\leq t\})_{0\leq t\leq T}$. Let $
\bF=(\cF_t)_{0\leq t\leq T}$ be the completed filtration of
$(\cF_t^0)_{0\leq t\leq T}$ with the $\mathbb{P}$-null sets of
${\cal F}$, hence $(\cF_t)_{0\leq t\leq T}$ satisfies the usual
conditions, $i.e.$, it is right continuous and complete.
Furthermore, let:

- ${\cal P}$ be the $\sigma$-algebra on $[0,T]\times \Omega$ of
$\bF$-progressively measurable sets ;

- ${\cal H}^{2,k}$ be the set of $\cal P$-measurable, $\R^k$-valued
processes $w=(w_t)_{t\leq T}$ such that
$\E[\int_0^T|w_s|^2ds]<\infty$ ;

- ${\cal S}^2$  be the set of $\cal P$-measurable, continuous,
$\R$-valued  processes ${w}=({w}_t)_{t\leq T}$ such that\\
$\E[\sup_{t\leq T}|{w}_t|^2]<\infty$.

\medskip

The problem of multiple switching can be described through an
example as follows. Assume  we have a plant which produces a
commodity, $e.g.$ a power station which produces electricity. Let
${\cal J} $ be the set of all possible activity modes of the
production of the commodity. A management strategy of the plant
consists, on the one hand, of the choice of a sequence of
nondecreasing stopping times $(\tau_n)_{n\geq1}$ (i.e. $\tau_n \leq
\tau_{n+1}$ and $\tau_0 = 0)$ where the manager decides to switch
the activity from its current mode to another one. On the other
hand, it consists of the choice of the mode $\xi_n$, a r.v. ${\cal
F}_{\tau_n}$-measurable with values in ${\cal J}$, to which the
production is switched at $\tau_n$ from its current mode. Therefore
the admissible management strategies of the plant are the pairs
$(\delta,\xi):=((\tau_n)_{n\geq 1},(\xi_n)_{n\geq 1})$ for which we
require also that $\P[\tau_n<T, \forall n\geq 0]=0$. This set is
called of {\it admissible strategies} and denoted by $\cal D$.

Next, assuming that the production activity is in mode 1 at the
initial time $t = 0$, let $(\a_t)_{t\leq T}$ denote the indicator of
the production activity's mode at time $t\in [0, T]$, i.e.,
\begin{equation}\label{proca}
\a_t=\ind_{[0,\tau_1]}(t)+\sum_{n\geq1}\xi_n
\ind_{(\tau_{n},\tau_{n+1}]}(t).
\end{equation}

Finally, for $i\in {\cal J}$, let $(\psi_i(t,\omega))_{t\leq T}$ be
a process of ${\cal H}^{2,1}$ which stands for the instantaneous
profit when the system is in state $i$ and for $i,j \in {\cal
J},\,i\neq j$, let $(g_{ij}(t,\omega))_{t\leq T}$ be a process of
${\cal S}^2$ which denotes the switching cost of the production at
time $t$ from current mode $i$ to another one $j$. If the plant is
run under the admissible strategy $(\delta,\xi)=((\tau_n)_{n\geq
1},(\xi_n)_{n\geq 1})$ the expected total profit is given by:
$$\begin{array}{l} J(\delta,\xi):=\E[\integ{0}{T}\psi_{\a_s}(s)ds
-\sum_{n\geq 1} g_{{\xi_{n-1}},{\xi_n}}(\tau_{n})
\ind_{[\tau_{n}<T]}].
\end{array}
$$
Therefore in several works authors are usually interested in either
finding an optimal strategy, $i.e$, a strategy $(\delta^*,\xi^*)$
such that $J(\delta^*,\xi^*)\ge J(\delta,\xi)$ for any
$(\delta,\xi)\in \cal D$ (\cite {Djehicheetal09, hutang}) or at
least in characterizing the quantity $\sup_{(\delta,\xi)\in {\cal
D}}J(\delta,\xi)$ in some specific cases (\cite{carmonaludkovski_1,
Hamadene_elasri}). This latter quantity is in a way the price of the
power plant in the energy market. \qed
\subsection{ Connection with systems of reflected BSDEs with oblique reflection}

In order to tackle the switching problem described above, we usually relate it to systems
 of reflected BSDEs with oblique reflection which we introduce below in the case we
need in order to deal also with the system of PDEs (\ref{sysvi1}).

Let $(t,x)\in [0,T]\times \R^k$ and $(X_s^{t,x})_{s\leq T}$ be the
solution of the following stochastic differential equation:
$$
dX_s^{t,x}=b(s,X_s^{t,x})ds+\sigma(s,X_s^{t,x})dB_s, \,\,s\in [t,T]
\mbox{ and }X_s^{t,x}=x \mbox{ for }s\in [0,t].$$ The solution of
this equation exists, is unique, since $b$ and $\sigma$ verify [H1],
and satisfies: \be\label{croix} \forall p\geq 1,\,\, \E[\sup_{s\leq
T}|X_s^{t,x}|^{p}]\leq C(1+|x|^p).\ee

Next let us introduce the solution of the system of reflected BSDEs
with oblique reflection associated with the deterministic functions
$((f_i)_{i\in {\cal J}},(g_{ij})_{i,j\in {\cal J}}, (h_i)_{i\in
{\cal J}})$ introduced in Section 2. The solution consists of $m$
triplets of processes $((Y^{i;t,x},Z^{i;t,x},K^{i;t,x}))_{i\in {\cal
J}}$, which is denoted, for convenience, by
$((Y^{i},Z^{i},K^{i}))_{i\in {\cal J}} $ and satisfies:   for any
$i\in {\cal J}$, \be\label{system}\left\{\begin{array}{l} Y^i, K^i
\in \cS^2, \,\,Z^i \in \cH^{2,d} \mbox{ and }K^i \mbox{
non-decreasing and }K^i_0=0 ;\\
Y^i_s=h_i(X^{t,x}_T)+\int_s^Tf_i(r,X^{t,x}_r,
Y^1_r,\dots,Y^m_r,Z^i_r)dr+K^i_T-K^i_s-\int_s^TZ^i_rdB_r,\,\,\forall\,\,s\leq
T\\
Y^i_s\geq \max_{j\in {\cal J}^{-i}}\{Y^j_s-g_{ij}(s,X^{t,x}_s)\},
\,\,\forall s\leq T\\
\int_0^T(Y^i_s-\max_{j\in {\cal
J}^{-i}}\{Y^j_s-g_{ij}(s,X^{t,x}_s)\})dK^i_s=0.\end{array}\right.\ee

We first provide an existence and uniqueness results of the
solution of (\ref{system}) and some of their properties as well.
\begin{theo}\label{hamjian} Assume that:

\noindent 1) the functions $(f_i)_{i\in \cJ}$ satisfy (H2)-(ii),
(iii) and (iv) ;

\noindent 2) For any $i,j\in {\cal J}$, the functions $g_{ij}$
(resp. $h_i$) verify (H3) (resp. (H4)).\ms

\no Then the system (\ref{system}) has a solution
$((Y^i,Z^i,K^i))_{i=1,m}$.
\end{theo}
\ms

\no \underline{Proof}: Since the above assumptions are not exactly
the same as the ones of Theorem 3.2 in \cite{hamzhang} we give
its main steps for sake of completeness. \ms

\no \underline{{\it Step 1}}: Let us consider the following BSDEs:
\be \label{solmax}\left\{\begin{array}{l}
\bar Y\in \cS^2, \bar Z\in \cH^{2,d}\\
\bar Y_s=\max_{i=1,m}h_i(X^{t,x}_T)+\int_s^T[\max_{i=1,m}
f_i](r,X^{t,x}_r,\bar Y_r,\dots, \bar Y_r,\bar Z_r)dr-\int_s^T\bar
Z_rdB_r,  s\leq T,
\end{array}
\right. \ee and \be \label{solmin} \left\{\begin{array}{l}
\underbar Y\in \cS^2, \underbar Z\in \cH^{2,d}\\
\underbar Y_s=\min_{i=1,m}h_i(X^{t,x}_T)+\int_s^T[\min_{i=1,m}
f_i](r,X^{t,x}_r,\underbar Y_r,\dots, \underbar Y_r,\underbar
Z_r)dr-\int_s^T\underbar Z_rdB_r, s\leq T.
\end{array}
\right. \ee Thanks to the result by Pardoux-Peng \cite{pardouxpeng}, the solutions of both (\ref{solmax}) and
(\ref{solmin}) exist. We next introduce the following sequences of
BSDEs defined recursively by: for any $i\in {\cal J}$,
$Y^{i,0}=\underbar Y$ and for $n\geq 1$ and $s\leq T$,
\be\label{sistemappro}\left\{\begin{array}{l} Y^{i,n},\, K^{i,n} \in
\cS^2, \,\,Z^{i,n} \in \cH^{2,d} \mbox{ and }K^{i,n} \mbox{
non-decreasing} ;\\
Y^{i,n}_s=h_i(X^{t,x}_T)+\int_s^Tf_i(r,X^{t,x}_r,
Y^{1,n-1}_r,\dots,Y^{i-1,n-1}_r, Y^{i,n}_r, Y^{i+1,n-1}_r,\dots,
Y^{m,n-1}_r,Z^{i,n}_r)dr\\\qquad\qquad\qquad+
K^{i,n}_T-K^{i,n}_s-\int_s^TZ^{i,n}_rdB_r;\\
Y^{i,n}_s\geq \max_{j\in {\cal
J}^{-i}}\{Y^{j,n-1}_s-g_{ij}(s,X^{t,x}_s)\};\\
\int_0^T(Y^{i,n}_s-\max_{j\in {\cal
J}^{-i}}\{Y^{j,n-1}_s-g_{ij}(s,X^{t,x}_s)\})dK^{i,n}_s=0.\end{array}\right.\ee
By an induction argument and the result by El-Karoui et al.
(\cite{Elkarouietal97}, Theorem 5.2), we claim that the processes
$(Y^{i,n},Z^{i,n},K^{i,n})$ exist for any $n\geq 1$. Next using the
comparison theorem of solutions of BSDEs (see e.g. Theorem 2.2 in
\cite{elkarouiquenezpeng}) we deduce that for any $\ij$,
$Y^{i,0}\leq Y^{i,1}.$ Now since $f_i$ satisfies the monotonicity
property (H2)-(iv) and using once more the comparison of solutions
of reflected BSDEs (see e.g. Theorem 4.1 in \cite{Elkarouietal97})
we obtain by induction that:
$$\forall\, n\geq 0 \mbox{ and }\ij,\,\,
Y^{i,n}\leq Y^{i,n+1}.$$ Since the processes $((\bar Y, \bar
Z,0))_{\ij}$ is a solution for the system of obliquely reflected
BSDEs associated with $(([\max_{i=1,m}
f_i](s,X^{t,x}_s,y^1,\dots,y^m,z))_{\ij},
(\max_{i=1,m}h_i(X^{t,x}_T))_{\ij},(g_{ij}(s,X^{t,x}_s))_{i,j\in
\cJ})$, an induction procedure and the repeated use of comparison theorem, which is justified taking into account that $f_i$ satisfies the
monotonicity property (H2)-(iv)), leads to
$$\forall n\geq 0, \; \; \forall \;i \in {\cal J},\,\,Y^{i,n}\leq \bar Y.$$
\ms

\no \underline{{\it Step 2}}: Using Peng's monotonic limit theorem (see Theorem
2.1 in \cite{Peng1999}), we deduce that for any $i\in {\cal J}$,
there exist: \ms

(i) a \cadlag (for right continuous with left limits) process $Y^i$
such that $Y^{i,n}\nearrow Y^i$ pointwisely ;

(ii) a process $Z^i$ of $\cH^{2,d}$ such that, at least for a
subsequence, $(Z^{i,n})_{n\geq 0}$ converges weakly to $Z^i$ in
$\cH^{2,d}$ and strongly in $L^p(dt\otimes dP)$ for any $p\in [1,2[$
;

(iii) a \cadlag non decreasing process $K^i$ such that for any
stopping time $\tau$, $(K^{i,n}_\tau)_{n\geq 0}$ converges to
$K^i_\tau$ in $L^p(dP)$ for any $p\in [1,2[$. \bs

\noindent Additionally in taking the limit in (\ref{sistemappro}),
the triple of processes $(Y^i,Z^i,K^i)$ satisfies: \be
\label{sistemappro1}\left\{\begin{array}{l}
Y^i_s=h_i(X^{t,x}_T)+\int_s^Tf_i(r,X^{t,x}_r,
Y^1_r,\dots,Y^m_r,Z^i_r)dr+K^i_T-K^i_s-\int_s^TZ^i_rdB_r\,\mbox{ and }\\
Y^i_s\geq \max_{j\in {\cal J}^{-i}}\{Y^j_s-g_{ij}(s,X^{t,x}_s)\},
\,\,\forall s\leq T.\end{array}\right.\ee Next let us consider the
following $m$ independent reflected BSDEs with a \cadlag barrier which, in addition, are
independent of each other (for existence and uniqueness results for such BSDEs, we refer to \cite{hamadenestochastics}): \be
\label{eqtilde}\left\{\begin{array}{l} \tilde
Y^i_s=h_i(X^{t,x}_T)+\int_s^Tf_i(r,X^{t,x}_r, Y^1_r,...,Y^{i-1}_r,
\tilde Y^i_r,Y^{i+1}_r, \dots, Y^m_r,Z^i_r)dr+\tilde K^i_T-\tilde
K^i_s-\int_s^T\tilde Z^i_rdB_r,\,\,s\leq
T\\
\tilde Y^i_s\geq \max_{j\in {\cal
J}^{-i}}\{Y^j_s-g_{ij}(s,X^{t,x}_s)\}, \,\,\forall s\leq
T\\
(\tilde Y^i_{s-}-\max_{j\in {\cal
J}^{-i}}\{Y^j_{s-}-g_{ij}(s,X^{t,x}_s)\} )d\tilde K^i_s=0,\,\,s\leq
T\end{array}\right.\ee where $\tilde Y^{i}_{t-}\,(\mb{resp.
}K^{i}_{t-})=\lim_{s\nearrow t}\tilde Y^i_s\, (\mb{resp. }
K^{i}_s)$, $t>0$. Henceforth in using once more comparison theorem
(e.g. \cite{{hamadenestochastics}}, Theorem 1.5) we have for any
$i\in {\cal J}$, $Y^{i,n}\leq \tilde Y^i$ and then $Y^{i}\leq \tilde
Y^i$. On the other hand using It\^o's formula with $((\tilde Y^i_t-
Y^i_t)^+)^2$, $t\leq T$, we obtain that $\tilde Y^i_t\leq Y^i_t$,
for any $t\leq T$. Actually this is true because $\tilde Y^i$ is the
minimal solution which satisfies system (\ref{sistemappro1}). It
implies that $\tilde Y^i=Y^i$, and then $Z^i=\tilde Z^i$, $\tilde
K^i=K^i$, for any $i\in {\cal J}$. Therefore $(Y^i,Z^i,K^i)$
satisfies not only the system (\ref{sistemappro1}) but also:
$$\label{contproofx}
\int_0^T(Y^i_{s-}-\max_{j\in {\cal
J}^{-i}}\{Y^j_{s-}-g_{ij}(s,X^{t,x}_s)\})dK^i_s=0, \,\,\forall i\in
\cJ.
$$
It remains to show that, for each $\ij$, $Y^i$ is continuous. For
this, let us assume that for some $i$ and $t_0$, $\D
Y^i_{t_0}:=(Y^i_{t_0}-Y^i_{t_0-})<0$, therefore $\D K^i_{t_0}>0$ and
then from (\ref{eqtilde}) we deduce the existence of $i_1\in {\cal
J}^{-i}$ such that
$Y^i_{{t_0}-}=Y^{i_1}_{{t_0}-}-g_{ii_1}({t_0},X^{t,x}_{t_0})$ and as
$Y^i_{{t_0}}\geq Y^{i_1}_{{t_0}}-g_{ii_1}({t_0},X^{t,x}_{t_0})$ then
$\D Y^{i_1}_{t_0}<0$. Repeating this reasoning as many times as
necessary we deduce the existence of a sequence of indices $(i_k)$
such that $i_k\in {\cal J}^{-i_{k-1}}$ and $$
Y^{i_{k-1}}_{{t_0}-}=Y^{i_k}_{{t_0}-}-g_{i_{k-1}i_k}({t_0},X^{t,x}_{t_0}).$$
As $\cJ$ is finite then there exists a sequence $j_1,..,j_{p-1}$
which are all different such that
$$g_{j_1j_2}({t_0},X^{t,x}_{t_0})+\dots +g_{j_{p-1}j_1}({t_0},X^{t,x}_{t_0})=0.$$
But this is contradictory with the assumption [H3]-(ii) and then
such a $t_0$ does not exist and $Y^i$ is continuous for any $i\in
\cJ$. Therefore $((Y^i,Z^i,K^i))_{i\in \cJ}$, is a solution for the
system of reflected BSDEs (\ref{system}).\qed \ms

We now give a remark related to comparison of the solutions of
system (\ref{system}) constructed in Theorem
\ref{hamjian}. Its proof is rather easy since an induction argument allows
to compare the solutions of the convergence schemes.
Actually we have:
\begin{remark} \label{comparison}: Let $(f'_i)_{\ij}$ (resp.$(g'_{ij})_{i,j\in \cJ}$, resp. $(h'_i)_{\ij}$)
be functions that satisfy (H2)-(ii), (iii), (iv) (resp. (H3), resp.
(H4)) and let  $((Y'^i,Z'^i,K'^i))_{\ij}$ be the solution of the
system of reflected BSDEs associated with
$((f'_i)_{\ij},(g'_{ij})_{i,j\in \cJ}, (h'_i)_{\ij})$ constructed as
in Theorem  \ref{hamjian}. If for any $i,j\in \cJ$ we have:
$$
f_i\leq f'_i,\,\, h_i\leq h'_i \mbox{ and }g_{ij}\geq g'_{ij}$$ then
for any $\ij$,
$$
Y^i\leq Y'^i.$$ In case of uniqueness of the solutions of those systems, this result reduces to the
comparison of the solutions.
\qed
\end{remark}

We now focus on the regularity properties of the solution of system
(\ref{system}) constructed in Theorem \ref{hamjian}.
\begin{prop}\label{sansz}Assume the assumptions of Theorem
\ref{hamjian} are fulfilled. Then there exist $lsc$ deterministic
functions $(v^i)_{\ij}$, defined on $\esp$, $\R$-valued and which
belong to $\Pi^g$, such that:
$$\forall \ij,\,\forall \,s\in [t,T],\,\, Y^i_s=v^i(s,X^{t,x}_s),$$where
$((Y^i,Z^i,K^i))_{\ij}$ is the solution of (\ref{system}) constructed in Theorem \ref{hamjian}.
\end{prop}
Proof: Actually under the hypotheses of Theorem \ref{hamjian}, there
exist deterministic continuous with polynomial growth functions
$\bar v(t,x)$ and $\underline{v}(t,x)$ with values in $\R$ such that
for any $s\in [t,T], \bar Y_s=\bar v(s,X^{t,x}_s)$ and $\underline{
Y}_s=\underline{v}(s,X^{t,x}_s)$ (\cite{elkarouiquenezpeng}, Th.4.1).

Next by induction and thanks to the result by El-Karoui et al.
(\cite{Elkarouietal97}, pp.729), there exist deterministic continuous
functions $v^{i,n}(t,x)$ in the class $\Pi^g$ such that for any $i\in \cJ$ and
$n\geq 0$, \be \label{karouiresult}
Y^{i,n}_s=v^{i,n}(s,X_s^{t,x}),\forall s\in [t,T],\ee
 the process $Y^{i,n}$ being defined as the unique solution of (\ref{sistemappro}) (see step 1, Theorem \ref{hamjian})
 As
$Y^{i,n}\leq Y^{i,n+1}\leq  \bar Y$ then for fixed $i$, the sequence
$(v^{i,n})_{n\geq 0}$ is non-decreasing and such that $v^{i,n}\leq \bar v$,
then it converges pointwisely to $v^i$. This latter function is therefore lower semi-continuous on $\esp$, of polynomial growth since
$\underline{v} \leq v^i\leq \bar v$ and finally for any $s\in [t,T],
Y^i_s=v^i(s,X^{t,x}_s).$ $\Box$ \ms

\begin{remark}\label{croiz}: Since $v^i$, $i\in \cJ$, belongs to
$\Pi^g$, then classically (see e.g. \cite{Elkarouietal97}) one can
show that for any $i\in \cJ$, $\|Z^{i}\|_{{\cal H}^{2,d}}(t,x)$ is
also of polynomial growth.\qed
\end{remark}

We next give a representation result for the solutions of system
(\ref{system}) and, as a by product, we obtain a uniqueness result
in some specific cases. Actually let fix
$\overrightarrow{u}:=(u^i)_{i=1,m}$ in $\cH^{2,m}$ and let us
consider the following system of reflected BSDEs with oblique
reflection:

$\forall i\in \cJ$, $\forall s\leq T$,
\be\label{eq3}\left\{\begin{array}{l} Y^{u,i},\,\, K^{u,i} \in
\cS^2, \,\,Z^{u,i} \in \cH^{2,d} \mbox{ and }K^{u,i} \mbox{
non-decreasing }(K^{u,i}_0=0) ;\\
Y^{u,i}_s=h_i(X^{t,x}_T)+\int_s^Tf_i(r,X^{t,x}_r,
\overrightarrow{u_r},Z^{u,i}_r)dr+K^{u,i}_T-K^{u,i}_s-\int_s^TZ^{u,i}_rdB_r;
\\
Y^{u,i}_s\geq \max_{j\in {\cal
J}^{-i}}\{Y^{u,j}_s-g_{ij}(s,X^{t,x}_s)\};\\
\int_0^T(Y^{u,i}_s-\max_{j\in {\cal
J}^{-i}}\{Y^{u,j}_s-g_{ij}(s,X^{t,x}_s)\})ds=0.\end{array}\right.\ee
Let $s\leq T$ be fixed and $i\in {\cal J}$ and let $\cD^i_s$ be the
following set of admissible strategies :
$$
\cD^i_s:=\{\a=((\tau_n)_{n\geq 0},(\xi_n)_{n\ge 0})\in \cD, \xi_0=i,
\t_0=0,\,\, \t_1\geq s \mbox{ and }E[(A^\a_T)^2]<\infty\}$$where for any $r\in [0,T]$,
$A^\a_r$ is the cumulative switching costs
up to time $r$, i.e.,
$$
A^\a_r:=\sum_{n\geq 1}g_{\xi_{n-1},\xi_n}(\t_n,X^{t,x}_{\tau_n})1_{[\t_n\leq r]}
\mbox{ for }r<T \mbox{ and } A_T=\lim_{r\rw T}A_r,\;\mathbb{P}\textrm{a.s}.
$$
Therefore and for any admissible strategy $\a$ we have:
$$
A^\a_T=\sum_{n\geq
1}g_{\xi_{n-1},\xi_n}(\t_n,X^{t,x}_{\tau_n})1_{[\t_n< T]}.
$$
Note that, by definition of the set $\cD$, $(\tau_{n}(\omega))$ is a stationnary sequence (for almost all $\omega$) and therefore
 the previous sum is finite, $\mathbb{P}$-almost surely.\\
Let us now consider a strategy $\a=((\tau_n)_{n\geq 0},(\xi_n)_{n\ge
0})\in \cD^i_s$ and let $(P^\a,Q^\a):=(P^\a_s,Q^\a_s)_{s\leq T}$ be
the solution of the following BSDE (which is not of standard type):
\be\label{eq4}\left\{\ba{l}
P^\a \mbox{ is RCLL and }E[\sup_{s\leq T}|P^\a_s|^2]<\infty, \,\,Q^\a\in \cH^{2,d};\\
P^\a_s=h_\a(X^{t,x}_T)+\int_s^Tf_\a(r,X^{t,x}_r,\overrightarrow{u_r},Q^\a_r)dr-\int_s^TQ^\a_rdB_r-(A^\a_T-A^\a_s),\,\forall
\,s\leq T,\ea\right. \ee with$$h_\a(x)=\sum_{n\geq
0}h_{\xi_{n}}(x)\ind_{[\tau_n<T\leq \tau_{n+1}[}\mbox{ and
}f_\a(t,x,(\zeta^i)_{\ij},q):=\sum_{n\geq
0}f_{\xi_{n}}(r,x,(\zeta^i)_{\ij},q)\ind_{[\tau_n\leq
r<\tau_{n+1}[}.$$ In setting up $\bar P^\a:=P^\a-A^\a$, we easily
deduce the existence and uniqueness of the process $(P^\a,Q^\a)$,
since $A^\a$ is adapted and $E[(A^a_T)^2]<\infty$, and the generator
as well as the terminal value of the transformed BSDE are standard.
\ms

We then have the following representation for the solution of
(\ref{eq3}) which is the main relationship between the value function of the optimal
switching problem and solutions of systems of reflected BSDEs with
oblique reflection. This result usually referred as the verification result is not new and has been already shown in
several contexts and under various assumptions.
\begin{theo}\label{chassegneux} Assume that for any $i,j\in\cJ$: \\
(i) $f_i$  satisfies (H2)-(ii), (iii) and for any fixed $(\overrightarrow{y},z)$ the mapping
$(t,x)\mapsto f_i(t,x,\overrightarrow{y},z)$ is $\cB (\esp)$-measurable ;\\
(ii) $g_{ij}$ (resp . $h_i$) satisfies (H3) (resp. (H4)).

\noindent Then the solution of system of BSDEs (\ref{eq3}) exists and satisfies:
\be\label{propprincip} \forall s\leq T, \,\,\forall
i\in \cJ,\,\,Y^{u,i}_s=\esssup_{\a\in \cD^i_s}P^\a_s.\ee
Thus the solution of (\ref{eq3}) is unique.
\end{theo}
\udl{Proof}: Thanks to Theorem \ref{hamjian} and
considering once more the same assumptions on the functions $(f_i,
g_{ij},h_i)_{i,j\in \cJ}$ the solution
$(Y^{u,i},Z^{u,i},K^{u,i})_{\ij}$ of system (\ref{eq3}) exists. Next using that $Y^{u,i}$ is solution
 to system (\ref{eq3}) and following the strategy $\a\in \cD_s^i$ in (\ref{system}), we obtain:
\be\label{secondterm} Y^{u,i}_s\geq
h_\a(X^{t,x}_T)+\int_s^Tf_\a(r,X^{t,x}_r,\overrightarrow{u_r},Z^\a_r)dr-
\int_s^TZ^\a_rdB_r-(A^\a_T-A^\a_s)+ \tilde K^\a_T\ee where
$$\tilde K^a_r=(K^i_{r}-K^i_{s}) \mbox{ if }t\leq r\leq \t_1,\,\tilde K^\a_r=\tilde
K^\a_{\t_n}+(K^{\xi_n}_{r}-K^{\xi_n}_{\t_n})\mbox{ if }\t_n<r\leq
\t_{n+1}\,;\, \forall r\leq T,\,\,Z^\a_r=\sum_{n\geq
0}Z^{\xi_i}1_{[\t_n\leq r<\t_{n+1}[}.$$ As $\tilde K^\a_T\geq 0$
then we have:
$$Y^{u,i}_s\geq P^\a_s,\,\,\forall\,\, \a\in \cD^i_s.$$
Note that the right-hand side in (\ref{secondterm}) is not a BSDE,
therefore we shall rather consider the equation satisfied by
$Y^{u,i} -P^{\alpha}$ where the pair $(P^\a,Q^\a)$ satisfies
(\ref{eq4}). Then using an equivalent change of probability we
deduce the previous inequality.

Next let $\a^*=(\t_n^*,\xi_n^*)_{n\geq 0}$ be the strategy defined
recursively as follows: $\t_0^*=0$, $\xi_0^*=i$ and for $n\ge 0$,
$$
\t^*_{n+1}=\inf\{s\geq \t_n^*, Y_s^{u,\xi_n^*}=\max_{j\in
\cJ^{-\xi_n^*}}(Y^{u,j}_s-g_{\xi_n^*,j}(s,\x_s))\}\wedge T$$ and
$$\xi_{n+1}^*=\mbox{argmax}_{j\in J^{-\xi_n^*}}\{Y^{u,j}_{\t^*_{n+1}}
-g_{\xi_n^*,j}(\t^*_{n+1},\x_{\t^*_{n+1}})\}.
$$
Let us show that $\a^*\in \cD^i_s$ and, for this, let us first prove that
$P[\t_n^*<T,\forall n\geq 0]=0.$ Actually assume the contrary i.e.
$P[\t_n^*<T,\forall n\geq 0]>0.$ Therefore thanks to definition of
$\t_n^*$, we have:
$$P[
Y_{\t_{n+1}^*}^{u,\xi_n^*}=Y^{u,\xi_{n+1}^*}_{\t_{n+1}^*}-g_{\xi_n^*,\xi_{n+1}^*}
({\t_{n+1}^*},\x_{\t_{n+1}^*}),\,\, \xi_n^*\in
\cJ^{-\xi_{n+1}^*},\,\,\forall n\geq 1]>0.$$ As $\cJ$ is  finite
then there is a state $i_0\in \cJ$ and a loop $i_0,i_1,...,i_k,i_0$
of elements of $\cJ$ such that $\mbox{card}\{i_0,i_1,...,i_k\}=k+1$
and
$$ P[
Y_{\t_{n+1}^*}^{u,i_l}=Y^{u,i_{l+1}}_{\t_{n+1}^*}-g_{i_l,i_{l+1}}
({\t_{n+1}^*},\x_{\t_{n+1}^*}),\,\,l=0,\dots,k,
\,(i_{k+1}=i_0),\,\forall n\geq 1]>0.$$ Therefore taking the limit
w.r.t. $n$ to obtain:
$$ P[
Y_{\t}^{u,i_l}=Y^{u,i_{l+1}}_{\t}-g_{i_l,i_{l+1}}
({\t},\x_{\t}),\,\,l=0,\dots,k, \,(i_{k+1}=i_0)]>0$$ where
$\t:=\lim_{n\rightarrow \infty}\t_n^*$. But this implies that
$$ P[
g_{i_0,i_1} ({\t},\x_{\t})+\dots+g_{i_k,i_0} ({\t},\x_{\t})=0]>0$$
which contradicts assumption $(H3)-(i)$. Thus we have
$P[\t_n^*<T,\forall n\geq 0]=0$.

Next it only remains to prove that $E[(A^{\a^*}_T)^2]<\infty$ and $\a^*$ is
optimal in $\cD^i_s$ for the switching problem (\ref{propprincip}).
Actually following the strategy $\a^*$ and since $Y^{u,i}$ solves the reflected BSDE (\ref{eq3}), it yields: for any $n\geq 1$,
\be\label{eq1} Y^{u,i}_s=
Y^{u,\xi_{n}^*}_{\t_n^*}+\int_s^{\t_n^*}f_{\a^*}(r,X^{t,x}_r,\overrightarrow{u_r},Z^{\a^*}_r)dr-
\int_s^{\t_n^*}Z^{\a^*}_rdB_r-A^{\a^*}_{\t_n^*}\ee noting that
$K^{\xi^*_n}_{r}-K^{\xi^*_n}_{\t_n}=0 \mbox{ holds for any $r$, }\t^*_n <r \leq
\t^*_{n+1}$. Taking now the limit w.r.t. $n$ in (\ref{eq1}) to
obtain: \be\label{eq2} Y^{u,i}_s=
h_{\a^*}(X^{t,x}_T)+\int_s^{T}f_{\a^*}(r,X^{t,x}_r,\overrightarrow{u_r},Z^{\a^*}_r)dr-
\int_s^{T}Z^{\a^*}_rdB_r-A^{\a^*}_{T}.\ee But taking into account
the assumptions (H4) and (H2)-(ii),(iii) satisfied by $h_i$ and
$f_i$ respectively and since $\overrightarrow{u}\in \cH^{2,m}$ and
$Z^{a^*}\in \cH^{2,d}$ and $(Y^{i})_{i}\in (\mathcal{S}^{2})^{m}$,
we deduce from (\ref{eq2}) that $E[(A^{\a^*}_{T})^2]<\infty$. It
follows that $\a^*\in \cD_s^i$ and $Y^{u,i}_s=P^{\a^*}_s$, thus
(\ref{propprincip}) holds and the solution of (\ref{eq3}) is unique.
\qed\ms

Next for $\orw{u}:=(u^i)_{i=1,m}\in {\cH}^{2,m}$ let us define by
\be\label{contraction}\Phi(\orw{u}):=(Y^{u,i})_{i=1,m}\ee where
$((Y^{u,i},Z^{u,i},K^{u,i}))_{i=1,m}$ is the solution of system
(\ref{eq3}) which exists and is unique if the assumptions of Theorem
\ref{chassegneux} are fulfilled. Note that when the processes
$(Y^{u,i})_{i=1,m}$ exist they belong to $({\cS}^2)^m$ and then
$\Phi$ is a mapping from ${\cH}^{2,m}$ to ${\cH}^{2,m}$.\ms

The following result, established by Chassagneux et al.
\cite{chassagneuxetal2010}, shows that $\Phi$ is a contraction
in ${\cH}^{2,m}$ when endowed with an appropriate equivalent norm.
Therefore the existence and uniqueness of a solution for
(\ref{system}) is deduced for general functions $f_i$ since, contrary to Theorem \ref{hamjian}, they are not supposed to satisfy any
monotonicity assumption. Actually we have:
\begin{theo} \label{theo3}Assume that for any $i,j\in \cJ$ the following assumptions are fullfiled:

(i) $f_i$ verifies (H2)-(ii),(iii) ;

(ii) $g_{ij}$ (resp. $h_i$) verifies (H3) (resp. (H4)).

\no Then we have: \ms

(a) For any $\orw{u}=(u^i)_{i=1,m},\,\orw{v}=(v^i)_{i=1,m} \in
{\cH}^{2,m}$,
\be\label{esti_theo3}\forall \, \ij, \,\forall s\leq T,\,
\E[|Y^{u,i}_s-Y^{v,i}_s|^2]\leq
C (\|\orw{u}-\orw{v}\|_{{\cH}^{2,m}}^2:=\E[\int_0^T\|\orw{u_r}-\orw{v_r}\|^2dr])\,;\ee

(b) The mapping $\Phi$ is a contraction when ${\cH}^{2,m}$ is
endowed with the following equivalent norm:
$$\|(u^i)_{i=1,m}\|_{\b_0}:=
\{\E[\int_0^Te^{\b_0 s}(\sum_{i=1,m}|u_s^i|^2)ds]\}^{\frac{1}{2}}, \,\,\mbox{ where
}\orw{u}=(u^i)_{i=1,m}\in {\cH}^{2,m}$$ for some appropriate $\b_0
\in \R$.
\end{theo}
\udl{Proof}: We provide it only for the sake of completeness since
it has been already given in \cite{chassagneuxetal2010}.  For $i\in
\cJ$, $\overrightarrow{u}$ and $\overrightarrow{v}\in \cH^{2,m}$ let
us set
$$ \varphi_i(r,X^{t,x}_r,z)= f_i(r,X^{t,x}_r, \overrightarrow{u_r},z)\vee
f_i(r,X^{t,x}_r, \overrightarrow{v_r},z),\,\, r\leq T,$$ and let us
consider the solution, denoted by $(\tilde Y^i,\tilde Z^i,\tilde
K^i)_{i\in \cJ}$, of the system of obliquely reflected BSDEs
associated with $(\varphi_i(r,X^{t,x}_r,z))_{i\in \cJ}$,
$(h_i)_{i\in \cJ}$ and $(g_{ij})_{i,j\in \cJ}$ which exists and is
unique by Theorem \ref{hamjian}. As shown in Theorem
\ref{chassegneux}, the following representation holds true:
$$\forall s\leq T,\,\,
\tilde Y^i_s=\esssup_{a\in \cD^i_s}\tilde P_s$$ where for $a\in
\cD^i_s$ the pair of processes $(\tilde P^a,\tilde Q^a)$ verifies:
$$\left\{\ba{l}
\tilde P^a \mbox{ is RCLL and }E[\sup_{\eta\leq T}|\tilde P^a_\et|^2]<\infty, \,\,\tilde Q^a\in \cH^{2,d};\\
\tilde P^a_\eta=h_a(X^{t,x}_T)+\int_\eta^T\varphi_a(r,\x_r,\tilde
Q^a_r)dr-\int_\eta^T\tilde Q^a_rdB_r-(A^a_T-A^a_\eta),\,\forall
\,\eta\leq T.\ea\right.
$$
Additionally an optimal strategy $\tilde a$ exists i.e. $\tilde
Y^i_s=\tilde P^{\tilde a}_s$. Note here that the dependence of
$P^{\tilde a}_s$ on $i$ is made through the strategy $\tilde a$ which
belongs to ${\cD}^i_s$. Now since for any $r\leq T$ and $z\in \R^d$,
$\varphi_i(r,X^{t,x}_r,z)\geq f_i(r,X^{t,x}_r,
\overrightarrow{u_r},z)$ and $\varphi_i(r,X^{t,x}_r,z)\geq
 f_i(r,X^{t,x}_r,
\overrightarrow{v_r},z)$ then by comparison and uniqueness (see
Remark \ref{comparison}) we have: \be \label{comp}Y^{u,i}\leq \tilde
Y^i\mbox{ and }Y^{v,i}\leq \tilde Y^i.\ee Next for $a\in\cD_s^i$,
let $(P^a_r,Q^a_r)_{r\leq T}$ be the solution of the non-standard
BSDE (\ref{eq4}) and let $(P'^a_r,Q'^a_r)_{r\leq T}$ be the solution
of the same non-standard BSDE with generator
$f_a(r,X^{t,x}_r,\orw{v_r},z)$. Then we have:
$$P^{\tilde a}_s\leq Y^{u,i}_s \leq \tilde
Y^i_s=\tilde P^{\tilde a}_s\mbox{ and }P'^{\tilde a}_s\leq
Y^{v,i}_s \leq \tilde Y^i_s=\tilde P^{\tilde a}_s$$ which implies,
\be \label{estimeyuv} |Y^{u,i}_s-Y^{v,i}_s|\leq |\tilde P^{\tilde
a}_s-P^{\tilde a}_s|+|\tilde P^{\tilde a}_s-P'^{\tilde a}_s|.
\ee But for any $\et \leq T$ we have:
$$\ba {l}\tilde P^{\tilde a}_\et-P^{\tilde a}_\et=\\
\int_\et^T\{f_{\tilde a}(r,X^{t,x}_r, \overrightarrow{u_r},\tilde
Q^{\tilde a}_r)\vee f_{\tilde a}(r,X^{t,x}_r,
\overrightarrow{v_r},\tilde Q^{\tilde a}_r)-f_{\tilde
a}(r,\x_r,\overrightarrow{u_r},Q^{\tilde
a}_r)\}dr-\int_\et^T\{\tilde Q^{\tilde a}_r-Q^{\tilde
a}_r\}dB_r,\ea
$$ and a similar equation is valid for $\tilde P^{\tilde
a}_\et-P'^{\tilde a}_\et.$
Next using It\^o's formula to obtain:
$$\ba {l}|\tilde P^{\tilde a}_\et-P^{\tilde a}_\et|^2+
\int_\et^T|\tilde Q^{\tilde a}_r-Q^{\tilde
a}_r|^2dr=-2\int_\et^T(\tilde P^{\tilde a}_r-P^{\tilde a}_r)\{\tilde Q^{\tilde a}_r-Q^{\tilde
a}_r\}dB_r\\\qquad\qquad +2
\int_\et^T(\tilde P^{\tilde a}_r-P^{\tilde a}_r)\{f_{\tilde a}(r,X^{t,x}_r, \overrightarrow{u_r},\tilde
Q^{\tilde a}_r)\vee f_{\tilde a}(r,X^{t,x}_r,
\overrightarrow{v_r},\tilde Q^{\tilde a}_r)-f_{\tilde
a}(r,\x_r,\overrightarrow{u_r},Q^{\tilde
a}_r)\}dr.\ea
$$
As for any $x,y\in \R$ we have $|x\vee y-y|\leq |x-y|$, then
$$\ba {l}|\tilde P^{\tilde a}_\et-P^{\tilde a}_\et|^2+
\int_\et^T|\tilde Q^{\tilde a}_r-Q^{\tilde
a}_r|^2dr=-2\int_\et^T(\tilde P^{\tilde a}_r-P^{\tilde a}_r)\{\tilde Q^{\tilde a}_r-Q^{\tilde
a}_r\}dB_r\\\qquad\qquad \qquad\qquad+2
\int_\et^T|\tilde P^{\tilde a}_r-P^{\tilde a}_r||f_{\tilde a}(r,X^{t,x}_r,
\overrightarrow{v_r},\tilde Q^{\tilde a}_r)-f_{\tilde
a}(r,\x_r,\overrightarrow{u_r},Q^{\tilde
a}_r)|dr.\ea
$$
Now classically we obtain the existence of a real constant $C\geq 0$
such that: \be\label{estimationu}\E[\sup_{\et\leq T}|\tilde
P^{\tilde a}_\et-P^{\tilde a}_\et|^2]\leq
C\E[\int_0^T\|\orw{u_r}-\orw{v_r}\|^2dr].\ee In the same way
considering $|\tilde P^{\tilde a}_\et-P'^{\tilde a}_\et|^2$ we
obtain a similar inequality as (\ref{estimationu}) where $P^{\tilde
a}$ is replaced by $P'^{\tilde a}$. Finally going back to
(\ref{estimeyuv}), squarring and taking the expectation, we obtain
the first estimate. \ms

Let us now show that $\Phi$ is a contraction. Let $\beta >0$ and let
us make use of It\^o's formula to obtain:
$$\ba {l}e^{\beta \et}|\tilde P^{\tilde a}_\et-P^{\tilde a}_\et|^2+
\int_\et^Te^{\beta r}|\tilde Q^{\tilde a}_r-Q^{\tilde
a}_r|^2dr=-2\int_\et^Te^{\beta r}(\tilde P^{\tilde a}_r-P^{\tilde
a}_r)(\tilde Q^{\tilde a}_r-Q^{\tilde a}_r)dB_r-{\beta}
\int_\et^Te^{\beta r}|\tilde P^{\tilde a}_r-P^{\tilde a}_r|^2dr
\\\qquad\qquad +2
\int_\et^Te^{\beta r}(\tilde P^{\tilde a}_r-P^{\tilde
a}_r)\{f_{\tilde a}(r,X^{t,x}_r, \overrightarrow{u_r},\tilde
Q^{\tilde a}_r)\vee f_{\tilde a}(r,X^{t,x}_r,
\overrightarrow{v_r},\tilde Q^{\tilde a}_r)-f_{\tilde
a}(r,\x_r,\overrightarrow{u_r},Q^{\tilde a}_r) \}dr,\,\,\eta \leq
T.\ea
$$
But $$| f_{\tilde a}(r,X^{t,x}_r, \overrightarrow{u_r},\tilde
Q^{\tilde a}_r)\vee f_{\tilde a}(r,X^{t,x}_r,
\overrightarrow{v_r},\tilde Q^{\tilde a}_r)-f_{\tilde
a}(r,\x_r,\overrightarrow{u_r},Q^{\tilde a}_r|\leq
C(|\overrightarrow{v_r}-\overrightarrow{u_r}|+|\tilde Q^{\tilde
a}_r-Q^{\tilde a}_r|)$$ where $C$ is the Lipschitz constant of $f$.
Therefore taking expectation in the previous equation and using both
inequalities $2Cxy\leq (Cx)^2+y^2$ and $2xy\leq
\frac{x^2}{\sqrt{\beta}}+\sqrt{\beta}y^2$ for any $x,y\in \R$, we
obtain: $\forall \,\et \leq T$,
$$
\E[e^{\beta \et}|\tilde P^{\tilde a}_\et-P^{\tilde a}_\et|^2]\leq
(C^2+C\sqrt{\beta}-\beta)\E[ \int_\et^Te^{\beta r}|\tilde P^{\tilde
a}_r-P^{\tilde
a}_r|^2dr]+\frac{C}{\sqrt{\beta}}\E[\int_\et^T\|\orw{u_r}-\orw{v_r}|^2dr]$$
But the same estimate can be obtained for $ \E[e^{\beta \et}|\tilde
P^{\tilde a}_\et-P'^{\tilde a}_\et|^2].$ Taking $\beta \geq
C^2+C\sqrt{\beta}$ and going back to (\ref{estimeyuv}) to obtain:
\be \label{sumi}\E[e^{\beta s}|Y^{u,i}_s-Y^{v,i}_s|^2]\leq
\frac{2C}{\sqrt{\beta}}\E[\int_0^T\|\orw{u_r}-\orw{v_r}|^2dr].\ee
Next summing for $i=1,m$ in (\ref{sumi}) and integrating w.r.t. $dt$
we obtain that:
$$
\|\Phi(\orw{u})-\Phi(\orw{v})\|^2_{\b}:= \E[\int_0^Te^{\beta
s}(\sum_{i=1,m}|Y^{u,i}_s-Y^{v,i}_s|^2)ds]\leq
\frac{2CT}{\sqrt{\beta}}\E[\int_0^T\|\orw{u_r}-\orw{v_r}\|^2dr].$$
Choosing now $\beta =\b_0 \geq 2\max\{(2CT)^2,C^2+C\sqrt{\beta}\}$
yields that $\Phi$ is a contraction in the Banach space
$({\cH}^{2,m},\|.\|_{\beta_0})$, therefore it has a fixed point
$(Y^i)_{i=1,m}$ which can be chosen continuous since
$\Phi((Y^i)_{i=1,m})\in (\cS^2)^m$. Thus the system of reflected
BSDEs with interconnected obstacles has a unique solution. \qed
\begin{remark}\label{normequiv}
Let $(Y^{i,0})_{\ij}$ be fixed processes of ${\cH}^{2,m}$ and for
$n\geq 1$ let us set $(Y^{i,n})_{\ij}=\Phi((Y^{i,n-1})_{\ij})$. Then
the sequence $((Y^{i,n})_{\ij})_{n\geq 0}$ converges in
$({\cH}^{2,m},\|.\|)$ to the unique solution of the system of
reflected BSDEs associated with $((f_i)_{i\in {\cal
J}},(g_{ij})_{i,j\in {\cal J}}, (h_i)_{i\in {\cal J}})$ since $\Phi$
is a contraction in $({\cH}^{2,m},\|.\|_{\beta_0})$ and the norms
$\|.\|_{\beta_0}$ and $\|.\|$ are equivalent.\qed
\end{remark}
\section{Uniqueness of the solution of the system of PDEs}
In this section we deal with the issue of uniqueness of the
solution of system (\ref{sysvi1}) and to do so, we first establish an auxiliary result which is a classical one in viscosity literature (see e.g. \cite{Pham2009}, pp. 76).
\begin{lemma} \label{modifsursolution}Let $(v_i(t,x))_{i=1,m}$ be a supersolution of the system
(\ref{sysvi1}), then for any $\gamma \geq 0$ there exists $\l_0
>0$ which does not depend on $\theta$ such that for any $\l\geq \l_0$ and $\theta >0$, the $m$-uplet
$(v_i(t,x)+\theta e^{-\l t}|x|^{2\g +2})_{i=1,m}$ is a supersolution
for (\ref{sysvi1}).
\end{lemma}
Proof: Without loss of generality we assume that the functions
$v_1,\dots,v_m$ are $lsc$. For sake of convenience, we do not use the previous
definition of a supersolution but an equivalent one (see e.g.
\cite{Crandalllions92}). Let $i\in \cJ$ be fixed and let
$\varphi \in \cC^{1,2}$ be such that the function
$\varphi-(v^i+\theta e^{-\l t}|x|^{2\g +2})$ has a local maximum in
$(t,x)$ which is equal to $0$. As $(v_i)_{i=1,m}$ is a supersolution
for (\ref{sysvi1}), then we have: $\forall \ij$,
$$\begin{array}{l} \min\left \{v_i(t,x)- \max\limits_{j\in{\cal
J}^{-i}}(-g_{ij}(t,x)+v_j(t,x))\right.;\\ \qquad -\partial_t
(\varphi (t,x)- \theta e^{-\l t}|x|^{2\g
+2})-\frac{1}{2}Tr(\sigma.\sigma
^\top(t,x)D^2_{xx}(\varphi (t,x)- \theta e^{-\l t}|x|^{2\g +2}))\\
\left.-b(t,x)^\top.D_x(\varphi (t,x)- \theta e^{-\l t}|x|^{2\g
+2})-f_i(t,x, (v^1,\dots,v^m)(t,x), \sigma^\top (t,x)D_x(\varphi (t,x)-
\theta e^{-\l t}|x|^{2\g +2}))\dis \right \}\geq 0\end{array}$$
which implies that \be \label{eqsursol}
\begin{array}{l}(v_i(t,x)+\theta e^{-\l t}|x|^{2\g +2})-
\max\limits_{j\in{\cal J}^{-i}}(-g_{ij}(t,x)+(v_j(t,x)+\theta e^{-\l
t}|x|^{2\g +2}))\\\qquad\qquad =v_i(t,x)- \max\limits_{j\in{\cal
J}^{-i}}(-g_{ij}(t,x)+v_j(t,x))\geq 0.\end{array} \ee On the other
hand:
$$\begin{array}{l}
-\partial_t (\varphi (t,x)- \theta e^{-\l t}|x|^{2\g
+2})-\frac{1}{2}Tr(\sigma.\sigma ^\top(t,x)D^2_{xx}(\varphi (t,x)-
\theta e^{-\l t}|x|^{2\g +2}))\\-D_x(\varphi (t,x)- \theta e^{-\l
t}|x|^{2\g +2}).b(t,x)-f_i(t,x, (v^1,\dots,v^m)(t,x),
\sigma^\top(t,x)D_x(\varphi (t,x)- \theta e^{-\l t}|x|^{2\g
+2}))\geq 0\end{array}
$$
and then \be\label{sursol}\begin{array}{l} -\partial_t \varphi
(t,x)-\frac{1}{2}Tr(\sigma.\sigma ^\top(t,x)D^2_{xx}\varphi
(t,x))-b(t,x)^\top.D_x\varphi (t,x)\\\qq\qq\qq\qq\qq-f_i(t,x,
(v^i(t,x)+\theta e^{-\l t}|x|^{2\g +2})_{i=1,m},
\sigma^\top(t,x)D_x\varphi (t,x))\\\qquad \geq\theta \l e^{-\l
t}|x|^{2\g +2}-\frac{1}{2}\theta e^{-\l t}Tr(\sigma.\sigma
^\top(t,x)D^2_{xx}|x|^{2\g +2})- \theta e^{-\l t} D_x(|x|^{2\g
+2}).b(t,x)\\\qquad \qquad +[f_i(t,x, (v^1,\dots,v^m)(t,x),
\sigma^\top(t,x)D_x(\varphi (t,x)- \theta e^{-\l t}|x|^{2\g
+2}))\\\qquad \qquad \qquad \qquad \qquad \qquad \qquad \qquad
-f_i(t,x, (v^1,\dots,v^m)(t,x), \sigma^\top(t,x)D_x\varphi (t,x))].
\end{array}
\ee But the last term in the right-hand side of this latter
inequality is equal to $$\theta e^{-\l
t}C^i_{t,x,\theta,\l}.\sigma^\top(t,x)D_x(|x|^{2\g +2}),$$ where
$C^i_{t,x,\theta,\l}$ is bounded by a constant independent of
$\theta$ since the function $f_i$ is uniformly Lipschitz w.r.t. $z$.
Therefore, taking into account the growth conditions on $b$ and
$\si$, there exists a constant $\l_0 \in \R^+$ which does not depend
on $\theta$ such that if $\l\geq \l_0$, the right-hand side of
(\ref{sursol}) is non-negative. Henceforth, noting that $i$ is
arbitrary in $\cJ$ together with (\ref{eqsursol}), we obtain that
$(v^i+\theta e^{-\l t}|x|^{2\g +2})_{i=1,m}$ is a viscosity
supersolution for (\ref{sysvi1}). \qed \bs
\\

We now establish the comparison property between supersolutions
and subsolutions of (\ref{sysvi1}) in the case when $f_i$ does
not depend on $(y^1,\dots,y^{i-1},y^{i+1},\dots,y^m)$ for any $i\in
\cJ$. Actually let us introduce the following assumption on the
functions $f_{i}$'s. \bs

\no {\bf [H5]}: For any $i\in \cJ$, the function $f_i$ does not
depend on $(y^1,\dots,y^{i-1},y^{i+1},\dots,y^m)$.
Note that this assumption replaces assumption (H2)(iv): this last one does not make sense
 any more when $f_i$ depends only on $y^i$, $z^i$.
\begin{prop}\label{comparisonviscosity} Assume both (H3) and (H4) and let suppose that the functions $f_i$, $i\in \cJ$, verify
(H2)-(i),(ii) and (H5). Let $(u^i(t,x))_{i=1,m}$ (resp.
$(v^i(t,x))_{i=1,m}$) be a subsolution (resp. a supersolution) of
the system (\ref{sysvi1}) which belongs to $\Pi^g$, then for any
$i\in \cJ$, we have:
$$
\forall (t,x)\in \esp, u^i(t,x)\leq v^i(t,x).$$
\end{prop}
Proof: First w.l.o.g we assume that $u_i$ (resp. $w^i$) is $usc$
(resp. $lsc$) for any $i\in \cJ$. Next let $\gamma>0$ and $C$ be such that
that for any $i\in \cJ$ we have:
$$
|u_i(t,x)|+|v^i(t,x)|\leq C(1+|x|^\g),\,\,\forall (t,x)\in \esp.$$
For sake of clarity, the proof is divided into two steps. \bs

\no \underline{Step 1}: To begin with we additionally assume that
the functions $f_i$, $i\in \cJ$, satisfy: \be\label{hypdec} \exists
\lambda <0\mbox{ s.t. }\forall t,x,z,\,\,\forall \,u\geq v,
\,\,f_i(t,x,u,z)-f_i(t,x,v,z)\leq \l (u-v).\ee

According to the previous lemma we know that for any $\theta
>0$ and $\l$ large enough $(v_i(t,x)+\theta e^{-\l t}|x|^{2\g +2})_{i=1,m}$ is also a
supersolution for (\ref{sysvi1}). Therefore it is enough to show
that for any $i\in \cJ$, we have:
$$
\forall (t,x)\in \esp, u^i(t,x)\leq v^i(t,x)+\theta e^{-\l
t}|x|^{2\g +2},$$ since in taking the limit as $\theta \rightarrow 0$
we obtain the desired result. So let us set
$w^{i,\theta,\l}(t,x)=v^i(t,x)+\theta e^{-\l t}|x|^{2\g +2}$,
$(t,x)\in \esp$ and we still denote $w^{i,\theta,\l}$ by $w^i$. Next assume there exists a point $(\bar
t,\bar x)\in \esp$ such that for $i\in \cJ$: $\max_{i\in
\cJ}(u^i(\bar t,\bar x)-w^i(\bar t,\bar x))>0.$ Next using
the growth condition there exists $R>0$ such that:
$$
\forall (t,x)\in \esp \mbox{ s.t. }|x|\geq R,\,\,u^i(
t,x)-w^i(t,x)<0.$$ Taking into account the values of the
subsolution and the supersolution at $T$, it implies that
\be\label{sgn}\begin{array}{l} 0<\max_{(t,x)\in \esp }\max_{i\in
\cJ}(u^i( t,x)-w^i(t,x))=\\\qquad \qq \max_{(t,x)\in [0,T[\times
B(0,R)}\max_{i\in \cJ}(u^i( t,x)-w^i(t,x))=\max_{i\in \cJ}(u^i(
t^*,x^*)-w^i(t^*,x^*)),\end{array} \ee where $B(0,R)$ is the open
ball in $\mathbb{R}^{k}$ centered in $0$ and of radius $R$ and $(t^*,x^*)\in [0,T[\times
B(0,R)$.

Now let us define $\tilde \cJ$ as: \be\label{tildej} \tilde
\cJ:=\{j\in \cJ, u_j(t^*,x^*)-w^j (t^*,x^*)=\max_{k\in \cJ}(
u_k(t^*,x^*)-w^k (t^*,x^*))\}. \ee First note that $\tilde \cJ$ is
not empty. Next for  $j\in \tilde \cJ$ and $n\geq 1$, let us define:
\be\label{eqphin}
\Phi^j_n(t,x,y):=u_j(t,x)-w^j(t,y)-\varphi_n(t,x,y), \,\,(t,x,y)\in
[0,T]\times \R^{2k}, \ee where:
$\varphi_n(t,x,y):=n|x-y|^{2\g+2}+|x-x^*|^{2}+|t-t^*|^2.$ Now let
$(t_n,x_n,y_n)\in [0,T]\times B'(0,R)^2$ be such that
$$
\Phi^j_n(t_n,x_n,y_n)=\max_{(t,x,y)\in [0,T]\times
B'(0,R)^2}\Phi^j_n(t,x,y),
$$which exists since $\Phi^j_n$ is $usc$ ($B'(0,R)$ is the closure of $B(0,R))$. Then we have:
\be
\label{estivisco}
\Phi^j_n(t^*,x^*,x^*)=u_j(t^*,x^*)-w^j(t^*,x^*)\leq
u_j(t^*,x^*)-w^j(t^*,x^*)+\varphi_n(t_n,x_n,y_n)\leq
u_j(t_n,x_n)-w^j(t_n,y_n).\ee The definition of $\varphi_n$ together with the growth condition of $u_j$ and
$w^j$ implies that $(x_n-y_n)_{n\geq 1}$ converges to $0$. Next
for any subsequence $((t_{n_l},x_{n_l},y_{n_l}))_l$  which converges
to $(\tilde t, \tilde x, \tilde x)$ we deduce from (\ref{estivisco})
that
$$
u_j(t^*,x^*)-w^j(t^*,x^*)\leq u_j(\tilde t,\tilde
x)-w^j(\tilde t,\tilde x),
$$ since $u_j$ is $usc$ and $w^j$ is $lsc$.
As the maximum of $u_j-w^j$ on $[0,T]\times B'(0,R)$ is reached in $(t^*,x^*)$ then
this last inequality is actually an equality. It implies, from the definition of $\varphi_n$ and
(\ref{estivisco}), that the sequence
$((t_{n},x_{n},y_{n}))_n$ converges to $(t^*,x^*,x^*)$ from which we deduce 
$$
n|x_n-y_n|^{2\g+2}\rightarrow_n 0 \mbox{ and
}(u_j(t_n,x_n),w^j(t_n,y_n))\rightarrow_n
(u_j(t^*,x^*),w^j(t^*,y^*)).
$$Actually this latter convergence holds since from (\ref{estivisco}) we first obtain,
$$
u_j(t^*,x^*)-w^j(t^*,x^*)\leq \liminf_n u_j(t_n,x_n)-\limsup_n
w^j(t_n,y_n),
$$
whereas the fact that $u_{j}$ (resp. $w_{j}$) is $usc$ (resp. $lsc$) gives
$$  \limsup_n u_j(t_n,x_n)-\liminf_n
w^j(t_n,y_n)\leq u_j(t^*,x^*)-w^j(t^*,x^*).$$
All these inequalities imply that
$$
\liminf_n u_j(t_n,x_n)=\limsup_n u_j(t_n,x_n)=u_j(t^*,x^*) \mbox{
and }\liminf_n w^j(t_n,x_n)=\limsup_n
w^j(t_n,x_n)=w^j(t^*,x^*).
$$
Next as in \cite{Ishiikoike91}, let us show by contradiction that
for some $k\in \tilde \cJ$ we have:
$$
u_k(t^*,x^*)>\max_{j\in \cJ^{-k}}(u_j(t^*,x^*)-g_{kj}(t^*,x^*)).$$
Actually suppose that for any $k\in \tilde \cJ$ we
have:
$$
u_k(t^*,x^*)\leq \max_{j\in
\cJ^{-k}}(u_j(t^*,x^*)-g_{kj}(t^*,x^*)),$$then there exists $j\in
\cJ^{-k}$ such that
$$
u_k(t^*,x^*)- u_j(t^*,x^*)\leq -g_{kj}(t^*,x^*).$$ But $w^k$ is a
supersolution of (\ref{sysvi1}), therefore we have
$$
w^k(t^*,x^*)\geq w^j(t^*,x^*)-g_{kj}(t^*,x^*)$$ and then
$$
u_k(t^*,x^*)- u_j(t^*,x^*)\leq -g_{kj}(t^*,x^*)\leq
w^k(t^*,x^*)-w^j(t^*,x^*).$$Therefore
$$
u_k(t^*,x^*)-w^k(t^*,x^*)=u_j(t^*,x^*)-w^j(t^*,x^*)$$
which implies that $j$ also belongs to $\tilde \cJ$ and
$$
u_k(t^*,x^*)- u_j(t^*,x^*)=-g_{kj}(t^*,x^*).$$ Repeating this
procedure as many times as necessary and since $\tilde \cJ$ is
finite we get the existence of a loop of indices $i_1,
...,i_p,i_{p+1}$ of $\tilde \cJ$ such that $i_1=i_{p+1}$ and
$$
g_{i_1,i_2}(t^*,x^*)+\dots+g_{i_p,i_{p+1}}(t^*,x^*)=0.$$ But this
contradicts the assumption (H3) on $g_{ij}$, $i,j\in \cJ$, whence
the desired result.$\Box$ \ms

To proceed let us consider $k\in \tilde \cJ$ such that:
\be\label{cil1} u_k(t^*,x^*)>\max_{j\in
\cJ^{-k}}(u_j(t^*,x^*)-g_{kj}(t^*,x^*)). \ee As the functions $u_j$,
$j\in \cJ$, are $usc$ and $g_{ij}$ are continuous, then there exists
$\rho>0$ such that for $(t,x)\in B((t^*,x^*),\rho)$ we have
$u_k(t,x)>\max_{j\in \cJ^{-k}}(u_j(t,x)-g_{kj}(t,x))$. Next and by construction it holds that $(t_n,x_n,u_k(t_n,x_n))_n\rw_n (t^*,x^*,u_k(t^*,x^*))$ and once
more since $u_j$ is $usc$ then for $n$ large enough we have:
\be\label{cil2}u_k(t_n,x_n)>\max_{j\in
\cJ^{-k}}(u_j(t_n,x_n)-g_{kj}(t_n,x_n)).\ee

Now applying Crandall-Ishii-Lions's Lemma (see e.g.
\cite{Crandalllions92} or \cite{Flemingetsoner}, pp.216) with
$\Phi_n^k$ (note that $k\in \tilde \cJ$ and (\ref{cil2}) is
satisfied) in $(t_n,x_n,y_n)$, there exist $(p^n_u,q^n_u,M^n_u)\in \bar
J^{2,+}u_k(t_n,x_n)$ and $( p^n_w,q^n_w,M^n_w)\in \bar
J^{2,-}w^k(t_n,y_n)$ such that:
$$p^n_u-p^n_w=\partial_t \varphi_n(t_n,x_n,y_n),\,\,
q^n_u \,\,(\mbox{resp. }q^n_w)\,=\partial_x
\varphi_n(t_n,x_n,y_n) \,(\mbox{resp. }=-\partial_y
\varphi_n(t_n,x_n,y_n)) \mbox{ and }
$$
\be \label{cil3}
\left (\begin{array}{ll} M_u^n&0\\
0&-N_w^n\end{array}\right )\leq A_n+\frac{1}{2n}A_n^2\ee where
$A_n=D^2_{(x,y)} \varphi_n(t_n,x_n,y_n)$. But
$$
\partial_t \varphi_nk(t,x,y)=2(t-t^*),\q
\partial_x
\varphi_n(t,x,y)=2(\gamma
+1)n(x-y)|x-y|^{2\gamma}+2(x-x^*)
\mbox{ and }
$$$$\partial_y
\varphi_n(t,x,y)= 2n(\gamma+1)(x-y)|x-y|^{2\gamma}.
$$
On the other hand,
$$D^2_{(x,y)} \varphi(t,x,y)=\left(\begin{array}{ll}
D^2_{xx} \varphi(t_n,x_n,y_n)&D^2_{xy}\tilde
\varphi(t,x,y)\\D^2_{xy}
\varphi(t,x,y)&D^2_{yy}
\varphi(t,x,y)\end{array}\right )
$$
with $$\begin{array}{l}D^2_{xx} \varphi(t,x,y)=
2n(\g+1)|x-y|^{2\g-2}\{|x-y|^2I_k+2\g(x-y)(x-y)^\top\}+2I_k,\end{array}
$$
$$\begin{array}{l}D^2_{yy} \varphi(t,x,y)=
2n(\g+1)|x-y|^{2\g-2}\{|x-y|^2I_k+2\g(x-y)(x-y)^\top\}\end{array}
$$
and finally
$$
\begin{array}{l}D^2_{xy} \varphi(t,x,y)=
-2n(\g+1)|x-y|^{2\g}I_k-4n\g(\g+1)(x-y)(x-y)^\top|x-y|^{2\g-2}.\end{array}
$$
As $(u_i)_{i\in \cJ}$ (resp. $(w^i)_{i\in \cJ}$) is a subsolution
(resp. supersolution) of (\ref{system}) and taking into account (\ref{cil2}), we obtain:
$$
-p^n_u-b(t_n,x_n)^\top.q^n_u-\frac{1}{2}
Tr[(\sigma\sigma^\top)(t_n,x_n)M^n_u]-f_j(t_n,x_n,u_j(t_n,x_n),\sigma
(t_n,x_n)^\top.q^n_u)\leq 0
$$
and
$$-p^n_w-b(t_n,y_n)^\top.q^n_w-\frac{1}{2}
Tr[(\sigma\sigma^\top)(t_n,y_n)M^n_w]-f_j(t_n,y_n,w^j(t_n,y_n),\sigma
(t_n,y_n)^\top.q^n_w)\geq 0.
$$
Making the difference between those two inequalities yields:
\be\label{cil5}\begin{array}{l}
-(p^n_u-p^n_w)-(b(t_n,x_n)^\top.q^n_u-b(t_n,y_n)^\top.q^n_w)-
\frac{1}{2}
Tr[\{\sigma\sigma^\top(t_n,x_n)M^n_u-\sigma\sigma^\top(t_n,y_n)M^n_w\}]\\\qquad\qquad-(f_j(t_n,x_n,u_j(t_n,x_n),\sigma^\top
(t_n,x_n).q^n_u)-f_j(t_n,y_n,w^j(t_n,y_n),\sigma^\top
(t_n,y_n).q^n_w))\leq 0.\end{array} \ee
But\be\label{estibn}\lim_{n\rw
\infty}\{|b(t_n,x_n)^\top.q^n_u-b(t_n,y_n)^\top.q^n_w)|+|\sigma(t_n,x_n)^\top.q^n_u-\sigma(t_n,y_n)^\top.q^n_w)|\}=0.
\ee As usual taking into account (\ref{cil3}) we have: \be
\label{estisn}
\limsup_nTr[\{\sigma\sigma^\top(t_n,x_n)M^n_u-\sigma\sigma^\top(t_n,y_n)M^n_w\}]\leq
0.\ee Finally taking the limit in (\ref{cil5}) and using the
assumption (\ref{hypdec}) to obtain:
$$
-\l (u_j(t^*,x^*)-w^j(t^*,x^*))\leq 0,
$$
which is contradictory with (\ref{sgn}) and then for any $j\in \cJ$
we have $u_j\leq w^j$. \qed \ms

\noindent \underline{Step 2}: The general case. \ms

\no Once more let $(u_j)_{j\in \cJ}$ (resp. $(w^j)_{j\in \cJ}$) be a
subsolution (resp. supersolution) of (\ref{sysvi1}). For $j\in \cJ$
let us set $\tilde u_j(t,x)=e^{\l t}u_j(t,x)$ and $\tilde
w^j(t,x)=e^{\l t}w^j(t,x)$, $(t,x)\in \esp$. Then $(\tilde
u_j)_{j\in \cJ}$ (resp. $(\tilde w^j)_{j\in \cJ}$) is a subsolution
(resp. supersolution) of the following system of variational
inequalities with oblique reflection: for any $i\in \cJ$,\be
\label{sysvi3} \left\{
\begin{array}{l}
\min\left \{\tilde v_i(t,x)- \max\limits_{j\in{\cal J}^{-i}}(-e^{\l
t}g_{ij}(t,x)+\tilde v_j(t,x))\right.,\\ \qquad \qquad
\left.-\partial_t\tilde v_i(t,x)- {\cal
L}\tilde v_i(t,x)+\l \tilde v_i(t,x)-e^{\l t}f_i(t,x,e^{-\l t}\tilde v^i(t,x),e^{-\l t}\sigma^\top(t,x).D_x \tilde v^i(t,x))\right\}=0;\\
\tilde v_i(T,x)=e^{\l T}h_i(x).
\end{array}\right.
\ee Actually let $i\in \cJ$ and let $\vp (t,x)$ be a
$C^{1,2}$-function such that $\vp - \tilde u_i$ has a minimum at
$(t,x)$ and $\vp(t,x)=\tilde u_i(t,x)$. Therefore $e^{-\l t}\vp-u_i$
has a minimum at $(t,x)$ and $e^{-\l t}\vp(t,x)=u_i(t,x)$. As $u_i$
is a subsolution then $\tilde u_i(T,x)\leq e^{\l T}h_i(x)$ and, \be
\label{sysvi2}
\begin{array}{l}
\min\left \{u_i(t,x)- \max\limits_{j\in{\cal
J}^{-i}}(-g_{ij}(t,x)+u_j(t,x))\right.,\\ \qquad \qquad
\left.-\partial_t(e^{-\l t}\vp(t,x))- {\cal L}(e^{-\l
t}\vp)(t,x)-f_i(t,x,e^{-\l t}\vp(t,x),\sigma^\top(t,x).D_xe^{-\l
t}\vp(t,x))\right\}\leq 0.
\end{array}
\ee Now if $$u_i(t,x)- \max\limits_{j\in{\cal
J}^{-i}}(-g_{ij}(t,x)+u_j(t,x))\leq 0,$$ then we have
$$
\tilde u_i(t,x)- \max\limits_{j\in{\cal J}^{-i}}(-e^{\l
t}g_{ij}(t,x)+\tilde u_j(t,x)),$$ and the viscosity subsolution
property is satisfied. If not, i.e., $u_i(t,x)-
\max\limits_{j\in{\cal J}^{-i}}(-g_{ij}(t,x)+u_j(t,x))>0$ then:
$$
-\partial_t(e^{-\l t}\vp(t,x))- {\cal L}(e^{-\l
t}\vp)(t,x)-f_i(t,x,e^{-\l t}\vp(t,x),\sigma^\top(t,x).D_xe^{-\l
t}\vp(t,x))\leq 0$$ which implies $$ -\partial_t\vp(t,x)+\l
\vp(t,x)- {\cal L}\vp(t,x)-e^{\l t}f_i(t,x,e^{-\l t}\tilde
u_i(t,x),e^{-\l t}\sigma^\top(t,x).D_x\vp(t,x))\leq 0.$$Therefore
once more the viscosity subsolution property is satisfied. As $i$ is
arbitrary in $\cJ$, then $(\tilde u_i)_{i\in \cJ}$ is a viscosity
subsolution for (\ref{sysvi2}).

In the same way one can show that $(\tilde w^j)_{j\in \cJ}$ is a
viscosity supersolution of (\ref{sysvi2}), whence the claim.
\ms

Next for $i\in \cJ$ let us set:$$F_i(t,x,u,z)=-\l u+e^{\l
t}f_i(t,x,e^{-\l t}u,e^{-\l t}z).$$Taking now $\l =
(1+\max_{i=1,m}C^i)$, where $C^i$ is the Lipschitz constant of $f_i$
w.r.t. to $u$,  to obtain that for $u\geq v$,
$
F_i(t,x,u,z)-F_i(t,x,v,z)\leq -(u-v).
$
It means that $F_i$ satisfies the assumption (\ref{hypdec}).
Therefore and according to the result proved in Step 1, for any $j\in
\cJ$, we have $\tilde u_j\leq \tilde w^j$ and also $u_j\leq w^j$.
The proof of the proposition is now complete. \qed \bs
\\

Next, thanks to Proposition \ref{comparisonviscosity}, we classically deduce both uniqueness and continuity results of
any solution of (\ref{sysvi1}) which belongs to $\Pi^g$. Actually if $(u_i)_{i\in \cJ}$ is
a solution then $(u^*_i)_{i\in \cJ}$ (resp. $(u_{i*})_{i\in \cJ}$)
is a subsolution (resp. supersolution)  for the system
(\ref{sysvi1}) in the class $\Pi^g$, then we deduce
that $u^*_{i}\leq u_{i*}$ and then $u^*_{i}=u_{i*}=u_i$, for any
$i\in \cJ$. Whence the continuity of $(u_i)_{i\in \cJ}$. \ms
To sum up, we have:
\begin{theo}\label{comparison2}Assume that Assumptions
(H2)-(i),(ii) and (H5) for $f_i$, $i\in \cJ$, are fulfilled. Then:

(i) The system of variational inequalities with inter-connected
obstacles (\ref{sysvi1}) has at most one solution in the class
$\Pi^g$ ;

(ii) If the solution in $\Pi^g$ exists, it is necessarily
continuous. \qed
\end{theo}

\section{Existence of a solution for the system}
\subsection{Case 1: $f_i$ depends only on $(y^i,z^i)$.}

In this specific case we have the following existence result:

\begin{theo}
Under (H2), (H3), (H4) and (H5), the following system of variational
inequalities with inter-connected obstacles \be \label{sysvi2x}
\left\{
\begin{array}{l}
\min\left \{v_i(t,x)- \max\limits_{j\in{\cal
J}^{-i}}(-g_{ij}(t,x)+v_j(t,x))\right.,\\ \qquad \qquad
\left.-\partial_tv_i(t,x)- {\cal
L}v_i(t,x)-f_i(t,x,v_i(t,x),\sigma^\top(t,x).D_xv^i(t,x))\right\}=0;\\
v_i(T,x)=h_i(x)
\end{array}\right.
\ee has a unique continuous solution $(v^i)_{i\in \cJ}$ in the class
$\Pi^g$.
\end{theo}
$Proof$: First note that the hypothesis (H2)-(iv) does not make any
sense in consideration with (H5). Now let $(v^i)_{i\in \cJ}$ be the
functions constructed in Prop. \ref{sansz} which are associated with
the solution of the system of reflected BSDEs with inter-connected
obstacles associated with $((f_i)_{\ij},(h_i)_{\ij},(g_{ij})_{i,j\in
\cJ})$, which both exist under (H2)-(H5). The functions $v^i$, $i\in
\cJ$, are of polynomial growth, thus locally bounded. Next
let us show that they are viscosity solutions for the system
(\ref{sysvi2x}).

For any $i\in \cJ$, $v^i$ is $lsc$, then $v^i=v^i_*$. So let us show
that the $m$-uplet $(v^i)_{i\in \cJ}$ is a viscosity supersolution
to (\ref{sysvi2x}). First note that for any $i\in \cJ$, \be
\label{eq6}v^i=\lim_{n\to \infty}\nearrow v^{i,n},\ee where
$v^{i,n}$, for $n\geq 1$, is defined in (\ref{karouiresult}). By
El-Karoui et al.'s result (\cite{Elkarouietal97}, Thm. 8.5),
$v^{i,n}$ is a viscosity solution of the following
variational inequality or PDE with obstacle: \be
\label{ineqvin}\left\{\ba{l} \min\{v^{i,n}(t,x)-\max_{j\in
\cJ^{-i}}[v^{j,n-1}(t,x)-g_{ij}(t,x)];\\\qquad\qquad\qquad-\partial_tv^{i,n}-\cL
v^{i,n}(t,x)-f_i(t,x,v^{i,n}(t,x),\sigma(t,x)^\top D_xv^{i,n}(t,x)\}=0;\\
v^{i,n}(T,x)=h_i(T,x). \ea\right. \ee Now let us fix $i\in \cJ$, let
$(t,x)\in [0,T[\times \R^k$ and $(p,q,M)\in \bar J^-v^i(t,x)$. By
(\ref{eq6}) and Lemma 6.1 in \cite{Crandalllions92}, there exist
sequences
$$
n_j\rightarrow \infty,\,\,\qquad (p_j,q_j,M_j)\in
J^-v^{i,n_j}(t_j,x_j)$$ such that:
$$
(t_j,x_j,v^{i,n_j}(t_j,x_j), p_j,q_j,M_j)\rightarrow
(t,x,v^{i}(t,x), p,q,M).$$Now from the viscosity supersolution
property for $v^{i,n_j}$ we have:
$$
-p_j-b(t_j,x_j)^\top q_j-\frac{1}{2}Tr(\si
\si^\top(t_j,x_j)M_j)-f_i(t_j,x_j,v^{i,n_j}(t,x),\sigma(t_j,x_j)^\top
q_j)\geq 0,
$$
and taking then the limit as $j\rightarrow \infty$ we obtain:
$$
-p-b(t,x)^\top q-\frac{1}{2}Tr(\si
\si^\top(t,x)M)-f_i(t,x,v^i(t,x),\sigma(t,x)^\top q)\geq 0.
$$
As $v^i(t,x)\geq \max_{j\in \cJ^{-i}}(v^j(t,x)-g_{ij}(t,x))$ and
$v^i(T,x)=h_i(x)$ then $v^i$ is a viscosity supersolution for the
following PDE with obstacle:
$$\left\{\ba{l}
\min\{v^{i}(t,x)-\max_{j\in
\cJ^{-i}}[v^{j}(t,x)-g_{ij}(t,x)];\\\qquad\qquad\qquad-\partial_tv^{i}-\cL
v^{i}(t,x)-f_i(t,x,v^i(t,x),\sigma(t,x)^\top D_xv^{i}(t,x)\}=0\\
v^{i}(T,x)=h_i(x). \ea\right. $$ Finally as $i$ is arbitrary in
$\cJ$ then the $m$-uplet $(v^1,\dots,v^m)$ is a viscosity
supersolution for the system of variational inequalities
(\ref{sysvi2x}). $\Box$ \ms

Next let us show that $(v^{i*})_{i\in \cJ}$ is a subsolution for
(\ref{sysvi2x}). First let us show that for any $i\in \cJ$,
$v^{i*}(T,x)=h_i(x)$. To begin with we are going to show that:
$$
\min\{v^{i*}(T,x)-h_i(x);v^{i*}(T,x)-\max_{j\in
\cJ^{-i}}(v^{j*}(T,x)-g_{ij}(T,x))\}=0.$$Actually
$$v^{i*}(T,x)=\limsup_{(t',x')\rightarrow (T,x),t'<T}v^i(t',x')\geq
\limsup_{(t',x')\rightarrow (T,x),t'<T}v^{i,n}(t',x'), \mbox{ for
any }n\geq 0$$ therefore \be \label{plusgrand}v^{i*}(T,x)\geq
v^{i,n}(T,x)=h_i(x)\ee since $v^{i,n}$ is continuous and at $t=T$ it
equals to $h_i(x)$. On the other hand for any $(t,x)$ we have:
$$v^i(t,x)\geq \max_{j\in
\cJ^{-i}}(v^{j}(t,x)-g_{ij}(t,x))\},
$$then
\be\label{plusgrand2} v^{i*}(T,x)\geq \max_{j\in
\cJ^{-i}}(v^{j*}(T,x)-g_{ij}(T,x)), \ee which with (\ref{plusgrand2})
imply that: \be
\label{posebruno}\min\{v^{i*}(T,x)-h_i(x);v^{i*}(T,x)-\max_{j\in
\cJ^{-i}}(v^{j*}(T,x)-g_{ij}(T,x))\}\geq 0.\ee

Let us now show that the left-hand side of (\ref{posebruno}) cannot
be positive. We first follow the same idea as in
\cite{brunoboucha}. So let us suppose that for some $x_0$, there is $\eps
>0$ such that:
$$
\min\{v^{i*}(T,x_0)-h_i(x_0);v^{i*}(T,x_0)-\max_{j\in
\cJ^{-i}}(v^{j*}(T,x_0)-g_{ij}(T,x_0))\}= 2\eps,$$ and let us
construct a contradiction. Let $(t_k,x_k)_{k\geq 1}$ be a sequence
in $\esp$ such that:
$$(t_k,x_k)\rightarrow (T,x_0) \mbox{ and
}v^{i}(t_k,x_k)\rightarrow v^{i*}(T,x_0) \mbox{ as }k\rightarrow
\infty.
$$
Since $v^{i,*}$ is $usc$ and of polynomial growth and taking into
account of (\ref{eq6}), we can find a sequence $(\varrho^n)_{n\geq
0}$ of functions of $\cC^{1,2}(\esp)$ such that $\vro^n\rightarrow
v^{i,*}$ and, on some neighbourhood $B_n$ of $(T,x_0)$ we have:
\be\label{eq7x}\min\{\vro^n(t,x)-h_i(x), \vro^n(t,x)-\max_{j\in
\cJ^{-i}}(v^{j*}(t,x)-g_{ij}(t,x))\}\geq \e,\,\,\forall (x,t)\in
B_n.\ee After possibly passing to a subsequence of $(t_k,x_k)_{k\geq
1}$ we can then assume that it holds on $B^n_k:=[t_k,T]\times
B(x_k,\d_n^k)$ for some $\d^n_k \in (0,1)$ small enough in such a
way that $B^n_k\subset B$. Now since $v^{i,*}$ is locally bounded
then there exists $\zeta >0$ such that $|v^{i,*}|\leq \zeta$ on
$B_n$. We can then assume that $\vro ^n\geq -2\zeta$ on $B_n$. Next
let us define:$$ \tilde \vro^n_k(t,x):=\vro^n(t,x)+\frac{4\zeta
|x-x_k|^2}{(\d^n_k)^2}+\sqrt{T-t}.$$ Note that $\tilde \vro^n_k\geq
\vro^n$ and \be\label{eq8}(v^{i,*}-\tilde \vro^n_k)(t,x)\leq -\zeta
\mbox{ for }(t,x)\in [t_k,T]\times
\partial B(x_k,\d^n_k).\ee
Next since $\partial_t(\sqrt{T-t})\rightarrow -\infty$ as
$t\rightarrow T$, we can choose $t_k$ large enough in front of
$\d^n_k$ and the derivatives of $\vro^n$ to ensure that
\be\label{eq9}-\cL \tilde \vro^n_k (t,x)\geq 0 \mbox{ on }B_n^k.\ee
Next let us consider the following stopping time
$\theta_n^k:=\inf\{s\geq t_k, (s,X^{t_k,x_k}_s)\in {B_n^k}^c\}\wedge
T$ where ${B_n^k}^c$ is the complement of ${B_n^k}$, and
$\vartheta_k:=\inf\{s\geq t_k, v^i(s,X^{t_k,x_k}_s)=\max_{j\in
\cJ^{-i}}(v^j(s,X^{t_k,x_k}_s)-g_{ij}(s,X^{t_k,x_k}_s))\}\wedge T.$
Applying now It\^o's formula to the process $(\tilde \vro^n_k(s, X_s))$ stopped at time $\theta_n^k \wedge \vartheta_k $ and taking into account (\ref{eq7x}),
(\ref{eq8}) and (\ref{eq9}) to obtain:
$$\ba{ll}
\tilde \vro^n_k(t_k,x_k)&=\E[\tilde \vro^n_k(\theta_n^k\wedge
\vartheta_k,X^{t_k,x_k}_{\theta_n^k \wedge
\vartheta_k})-\int_{t_k}^{\theta_n^k\wedge \vartheta_k}
\cL \tilde \vro^n_k (r,X_r^{t_k,x_k})dr]\\
{}& \geq \E[\tilde
\vro^n_k(\theta_n^k,X^{t_k,x_k}_{\theta_n^k})\ind_{[\theta_n^k\leq
\vartheta_k]}
+\tilde\vro^n_k(\vartheta_k,X^{t_k,x_k}_{\vartheta_k})\ind_{[\vartheta_k<
\theta_n^k]}]\\{}&= \E[\{\tilde
\vro^n_k(\theta_n^k,X^{t_k,x_k}_{\theta_n^k})\ind_{[\theta_n^k<T]}+
\tilde
\vro^n_k(T,X^{t_k,x_k}_{T})\ind_{[\theta_n^k=T]}\}\ind_{[\theta_n^k\leq
\vartheta_k]}
+\tilde\vro^n_k(\vartheta_k,X^{t_k,x_k}_{\vartheta_k})\ind_{[\vartheta_k<
\theta_n^k]}]
\\
{}&\geq
\E[\{(v^{i*}(\theta_n^k,X^{t_k,x_k}_{\theta_n^k})+\zeta)\ind_{[\theta_n^k<T]}+
(\eps
+h_i(T,X^{t_k,x_k}_{T}))\ind_{[\theta_n^k=T]}\}\ind_{[\theta_n^k\leq
\vartheta_k]} \\{}&\qquad \qquad +\{\e+\max_{j\in
\cJ^{-i}}(v^{j*}(\vartheta_k,X^{t_k,x_k}_{\vartheta_k})-
g_{ij}(\vartheta_k,X^{t_k,x_k}_{\vartheta_k}))\}\ind_{[\vartheta_k<
\theta_n^k]}]\\
{}&\geq
\E[\{(v^{i}(\theta_n^k,X^{t_k,x_k}_{\theta_n^k})+\zeta)\ind_{[\theta_n^k<T]}+
(\eps
+h_i(T,X^{t_k,x_k}_{T}))\ind_{[\theta_n^k=T]}\}\ind_{[\theta_n^k\leq
\vartheta_k]} \\{}&\qquad \qquad +\{\e+\max_{j\in
\cJ^{-i}}(v^{j}(\vartheta_k,X^{t_k,x_k}_{\vartheta_k})-
g_{ij}(\vartheta_k,X^{t_k,x_k}_{\vartheta_k}))\}\ind_{[\vartheta_k<
\theta_n^k]}] \mbox{ since }v^{i*}\geq v^i\\
{}&\geq \E[v^{i}(\theta_n^k\wedge
\vartheta_k,X^{t_k,x_k}_{\theta_n^k \wedge
\vartheta_k})]+\zeta\wedge \eps\\{}&
=\E[v^i(t_k,x_k)-\int_{t_k}^{\theta_n^k\wedge
\vartheta_k}f_i(s,X_s^{t_k,x_k},v^i(s,X_s^{t_k,x_k}),Z^{i,t_k,x_k}_s)ds]+\zeta\wedge
\eps \ea$$ since on $[t_k,\vartheta_k]$, $dK^{i,t,x}=0$. Finally
since $v^i$ and $\|Z^{i}\|_{{\cal H}^{2,d}}(t,x)$ belong to $\Pi^g$ (see
Prop.\ref{sansz} and Remark \ref{croiz}) and taking into account
(\ref{croix}) and assumption (H2)-(iii), we easily deduce that$$\lim_{k\rightarrow \infty}\E[\int_{t_k}^{\theta_n^k\wedge
\vartheta_k}f_i(s,X_s^{t_k,x_k},v^i(s,X_s^{t_k,x_k}),
Z^{i,t_k,x_k}_s)ds]=0 .$$ Therefore taking the limit in the previous
inequalities yields:
$$\label{eq10} \lim_{k\rightarrow
\infty}\tilde \vro^n_k(t_k,x_k)=\lim_{k\rightarrow
\infty}\vro^n(t_k,x_k)+\sqrt{T-t_k}=\vro^n(T,x_0)\geq
\lim_{k\rightarrow \infty} v^{i}(t_k,x_k)+\zeta\wedge
\eps=v^{*i}(T,x_0)+\zeta\wedge \eps.$$ But this is a contradiction since $\vro^n\rightarrow v^{i*}$ pointwisely as $n\rightarrow
\infty$. Thus for any $(t,x)\in \esp$ we have:
$$
\min\{v^{i*}(T,x)-h_i(x);v^{i*}(T,x)-\max_{j\in
\cJ^{-i}}(v^{j*}(T,x)-g_{ij}(T,x))\}=0.$$

Let us now show that $v^{i*}(T,x)=h_i(x)$. So suppose that
$v^{i*}(T,x)>h_i(x)$, then for for some $j\in \cJ^{-i}$ we have
$$
v^{i*}(T,x)=v^{j*}(T,x)-g_{ij}(T,x).$$ But once more we have
$v^{j*}(T,x)>h_j(T,x)$. Otherwise, i.e. if $v^{j*}(T,x)=h_j(T,x)$,
we would have:
$$
h_i(x)<h_j(x)-g_{ij}(T,x),$$ which is contradictory with (H4).
Therefore there exists $\ell\in \cJ^{-j}$ such that:
$$
v^{j*}(T,x)=v^{\ell*}(T,x)-g_{j\ell}(T,x) \mbox{ and then
}v^{i*}(T,x)=v^{\ell*}(T,x)-g_{j\ell}(T,x)-g_{ij}(T,x).$$ Repeating
this reasoning as many times as necessary we obtain a sequence of
different indices $i_1,...,i_l$ such that $$
v^{i_1*}(T,x)=v^{i_1*}(T,x)-(g_{i_1i_2}(T,x)+\dots+g_{i_{l-1}i_l}(T,x)+g_{i_{l}i_1}(T,x)),$$
which is contradictory with (H3)-(ii). Thus for any $i\in \cJ$ we
have:
$$\forall x\in \R^k,\,\,
v^{i*}(T,x)=h_i(x).\qquad \Box$$

Let us now show that $(v^{i*})_{i\in \cJ}$ is a subsolution to
(\ref{sysvi2x}). First note that since $v^{i,n}\nearrow v^i$ and
$v^{i,n}$ is continuous then we have (see e.g. \cite{Pham2009},
pp.91) \be\label{def10} v^{i*}(t,x)=\lim_{n\rightarrow \infty}{
\sup}^* v^{i, n}(t,x)=\lim \sup_{n\rightarrow \infty, t'\rw t,x'\rw
x} v^{i, n}(t',x').\ee Next let us fix $i\in \cJ$ and let $(t,x)\in
[0,T[\times \R^k$ be such that
\be\label{eq11}v^{i*}(t,x)-\max_{\ell\in
\cJ^{-i}}(v^{\ell*}(t,x)-g_{i\ell}(t,x))>0.\ee Let $(p,q,M)\in \bar
J^+v^{i*}(t,x)$. By (\ref{def10}) and Lemma 6.1 in
\cite{Crandalllions92}, there exist sequences
$$
n_j\rightarrow \infty,\,\,\qquad (p_j,q_j,M_j)\in
J^+v^{i,n_j}(t_j,x_j)$$ such that:
$$
\lim_{j\rw \infty}\,\,(t_j,x_j,v^{i,n_j}(t_j,x_j),
p_j,q_j,M_j)=(t,x,v^{i*}(t,x), p,q,M).$$ Now from the viscosity
subsolution property for $v^{i,n_j}$ at $(t_j,x_j)$ (see
\ref{ineqvin}), for any $j\geq 0$, we have: \be\label{ssvn}\ba{l}
\min\{v^{i,n_j}(t_j,x_j)-\max_{\ell\in
\cJ^{-i}}(v^{\ell,n_j-1}(t_j,x_j)-g_{i\ell}(t_j,x_j));\\\qquad\qquad
-p_j-b(t_j,x_j)^\top q_j-\frac{1}{2}Tr(\si
\si^\top(t_j,x_j)M_j)-f_i(t_j,x_j,v^{i,n_j}(t_j,x_j),\sigma(t_j,x_j)^\top
q_j)\}\leq 0.\ea \ee Next the definition of $v^{\ell*}$ implies that
$$v^{\ell*}(t,x)\geq \lim\sup_{j\rw \fin}v^{\ell,n_j}(t_j,x_j),
$$
therefore by (\ref{eq11}), there exists $j_0\geq 0$, such that if
$j\geq j_0$ we have
$$v^{i,n_j}(t_j,x_j)>\max_{\ell \in
\cJ^{-i}}(v^{\ell,n_j}(t_j,x_j)-g_{ij}(t_j,x_j)).$$ Then
(\ref{ssvn}) implies that, for any $j\geq j_0$,
$$
-p_j-b(t_j,x_j)^\top q_j-\frac{1}{2}Tr(\si
\si^\top(t_j,x_j)M_j)-f_i(t_j,x_j,v^{i,n_j}(t_j,x_j),\sigma(t_j,x_j)^\top
q_j)\le 0.
$$
Taking the limit as $j\rightarrow \fin$ we deduce that
$$
-p-b(t,x)^\top q-\frac{1}{2}Tr(\si \si^\top
(t,x)M)-f_i(t,x,v^{i*}(t,x),\sigma(t,x)^\top q)\le 0
$$
since $f_i$ is uniformly continuous in $(t,x)$ and Lipschitz in
$(y^i,z)$. Then
$$\ba{l}
\min\{v^{i*}(t,x)-\max_{\ell\in
\cJ^{-i}}(v^{\ell*}(t,x)-g_{i\ell}(t,x));\\\qquad\qquad
-p-b(t,x)^\top q-\frac{1}{2}Tr(\si \si^\top
(t,x)M)-f_i(t,x,v^{i*}(t,x),\sigma(t,x)^\top q)\}\leq 0\ea
$$
which means that $v^{i*}$ is a viscosity subsolution for
(\ref{sysvi2x}). Thus the $m$-uplet $(v^i)_{i\in \cJ}$ is a solution
for (\ref{sysvi2x}) and Theorem \ref{comparison2} implies that it is
continuous and unique. \qed

As a by-product we obtain:
\begin{corollary}\label{representation}
Under (H1)-(H5), there exist deterministic continuous functions
$(v^i(t,x))_{\ij}$ which belong to $\Pi^g$ unique solution of
(\ref{sysvi2}) and such that the unique solution of the system of
reflected BSDEs with inter-connected obstacles associated with
$((f_i)_{\ij},(h_i)_{\ij},(g_{ij})_{i,j\in \cJ})$ has the following
representation: \be\label{caract} \forall \ij, \forall (t,x)\in
\esp, \forall s\in [t,T], \quad
Y^{i;t,x}_s=v^i(s,X^{t,x}_s).\,\,\,\qed\ee
\end{corollary}

\subsection{Case 2: Existence and uniqueness in
the general setting}

In this section we deal with the issue of existence and uniqueness
of the solution of (\ref{sysvi1}) in its general form, i.e., when
the functions $f_i$, $i\in \cJ$, depend not only on $y^i$ but also
on the other components $(y_1,\dots,y_m)$. More precisely, we
successively consider the cases when either $f_{i}$ or $-f_{i}$
satisfies (H2)-(iv) for any $i\in \cJ$. Considering the first case,
we have:
\begin{theo}\label{exist_increasingcase} Under (H2), (H3), (H4), the system of variational
inequalities with inter-connected obstacles (\ref{sysvi1}) has a
continuous solution $(v^1,\dots,v^m)$ in the class $\Pi^g$.
\end{theo}
Proof: We first prove existence. Let $\l \in \R$ and for $i\in\cJ$,
let $F_i$ be the function defined by: \be
\label{expgrdf}F_i(t,x,y^1,\dots,y^m,z)=e^{\l t}f_i(t,x,e^{-\l
t}y_1,...,e^{-\l t}y_m,e^{-\l t}z_i)-\l y_i. \ee Since $f_i$ is
uniformly Lipschitz w.r.t. $y_i$ then $F_i$ is so and for $\l$ small
enough ($\l<0$) the function $F_i$ is non-decreasing in all
variables $(y_1,...,y_m)$. Next thanks to Theorems \ref{hamjian} and
\ref{chassegneux}, there exist processes $(Y^i,Z^i,K^i)_{\ij}$
solution of the system of reflected BSDEs with interconnected
obstacles associated with $((F_i(t,x,y^1,...,y^m,z))_{i\in \cJ},
(e^{\l T}h_i(X^{t,x}_T))_{i\in \cJ}, (e^{\l t}g_{ij})_{i,j\in
\cJ})$. Additionally thanks to Proposition \ref{sansz}, there exist
deterministic $lsc$ functions $v^i$, $\ij$, such that:
\be\label{eqyschcr}\fl \ij,\,\forall t\leq T, \forall s\in [t,T],
Y^i_s=v^i(s,X^{t,x}_s).\ee

Let us now analyze the decreasing scheme. First let us consider $(\bar Y,\bar Z)$ the solution of the following standard BSDE:
\be
\label{bary}\left\{\ba{l}
\bar Y\in \cS^2, \,\,\bar Z\in \cH^{2,d}\\
\bar Y_s=\max_{i=1,m}e^{\l T}h_i(X^{t,x}_T)+\int_s^\top
[\max_{i=1,m}F_i](r,X^{t,x}_r,\bar Y_r,\dots,\bar Y_r,\bar
Z_r)dr-\int_s^\top\bar Z_rdB_r\,\,s\leq T. \ea\right.\ee Next for
any $\ij$, let us set $Y^{i,0}=\bar Y$ and for $n\geq 1$ let us
define $(Y^{i,n},Z^{i,n},K^{i,n})$ by:
\be\label{sistemappro3}\left\{\begin{array}{l} Y^{i,n},\, K^{i,n}
\in \cS^2, \,\,Z^{i,n} \in \cH^{2,d} \mbox{ and }K^{i,n} \mbox{
non-decreasing} ;\\
Y^{i,n}_s=e^{\l T}h_i(X^{t,x}_T)+\int_s^T F_i(r,X^{t,x}_r,
(Y^{i,n-1}_r)_{i=1,m},Z^{i,n}_r)dr+
K^{i,n}_T-K^{i,n}_s-\int_s^TZ^{i,n}_rdB_r;\\
Y^{i,n}_s\geq \max_{j\in {\cal J}^{-i}}\{Y^{j,n}_s-e^{\l
s}g_{ij}(s,X^{t,x}_s)\};\\
\int_0^T(Y^{i,n}_s-\max_{j\in {\cal J}^{-i}}\{Y^{j,n}_s-e^{\l
s}g_{ij}(s,X^{t,x}_s)\})dK^{i,n}_s=0.\end{array}\right.\ee
First note that the existence of $\bar Y$ is obvious by
Pardoux-Peng's result \cite{pardouxpeng} and then we easily deduce by
induction:

(i) for any $n\geq 1$, there exists a unique $m$-uplet of processes
$(Y^{i,n},Z^{i,n},K^{i,n})$, $\ij$, ;

(ii) for any $n\geq 1$, there exist deterministic continuous functions $v^{i,n}$,
$\ij$, such that: \be\label{schdecr}\fl \ij,\,\forall t\leq T,
\forall s\in [t,T], Y^{i,n}_s=v^{i,n}(s,X^{t,x}_s);\ee

(iii) for any $\ij$ and $n\geq 0$, $Y^{i,n+1}\leq Y^{i,n}$ and
$v^{i,n}\geq v^{i,n+1}$ ;

(iv) for any $\ij$, the sequence $(Y^{i,n})_n$ converges to $Y^i$ in
${\cH}^{2,1}$, where such a process $Y^{i}$ has been introduced in equation (\ref{eqyschcr}).\ms

\no Actually for (i), we just need to use Theorem \ref{hamjian} and \ref{chassegneux} since
$Y^{i,0}$ exists. To show (ii), we use induction: obviously the
property is valid for $n=0$ (see e.g. \cite{elkarouiquenezpeng}, Thm.4.1). Next if the property
holds for some $n$ then it holds also for $n+1$ in using Corollary
\ref{representation}. As for (iii), the property is true for $n=0$
since $((\bar Y^i,\bar Z^i,\bar K^i)=(\bar Y, \bar Z,0))_{\ij}$ is
the unique solution of the system
associated with \\
$(\bar f_i:=[\max_{i=1,m}F_i](r,X^{t,x}_r,y^1,...,y^m,z))_{i\in
\cJ}, (\bar h_i:=\max_{i=1,m}e^{\l T}h_i(X^{t,x}_T))_{i\in \cJ},
(e^{\l t}g_{ij})_{i,j\in \cJ})$ and then it is just enough to use
the comparison result of Remark \ref{comparison}. Next if the
property is valid for some $n$ then it is also valid for $n+1$, in
using once more comparison since $F_i(t,x,(y^i)_{\ij},z)$, $\ij$, is
non-decreasing in $(y^i)_{\ij}$. Finally (iv) holds true because the
mapping $\Phi$ defined in (\ref{contraction}) is a contraction in
$\cH^{2,m}$ and we obviously have:
$$\forall n\geq 0, (Y^{i,n+1})_{i=1,m}=\Phi((Y^{i,n})_{i=1,m}).$$

Let us now show that the deterministic functions $(v^i)_{\ij}$ of
(\ref{eqyschcr}) are continuous. Recall that thanks to Proposition \ref{sansz} each $ v^{i}$ is the limit of some increasing sequence and hence it is $lsc$: therefore,
 it is enough to show that they are upper semicontinuous. But from the estimate given in (\ref{esti_theo3}), we
know that: $$ \forall (t,x)\in
\esp,\,\,|v^i(t,x)-v^{i,n}(t,x)|=\E[|Y^i_t-Y^{i,n}_t|]\leq C
\|(Y^i)_{\ij}-(Y^{i,n-1})_{\ij}\|_{{\cH}^{2,m}}.$$ As $\lim_{n\rw
\infty}\|(Y^i)_{\ij}-(Y^{i,n-1})_{\ij}\|_{{\cH}^{2,m}}=0$ then for
any $\ij$, the sequence $(v^{i,n})_n$ defined in (\ref{sistemappro3}) converges pointwisely and decreasingly (from (\ref{sistemappro3})-(iii)) to
$v^i$. As a decreasing limit of continuous functions, $v^i$ is $usc$ and then continuous. It follows that the
solution $(y^i,z^i,k^i)_{i\in \cJ}$ of the system of reflected BSDEs
with oblique associated with\\
$((f_i(s,X^{t,x}_s,y^1,\dots,y^m,z))_{i\in
\cJ},(h_i(X^{t,x}_T))_{i\in \cJ}, (g_{ij}(s,X^{t,x}_s))_{i,j\in
\cJ})$ has the following representation:
$$
\forall \ij, \forall t\leq T, \forall s\in [t,T],\forall x\in
\R^k,\,\,\,y^i_s=e^{-\l s}v^i(s, X^{t,x}_s).$$ As $v^i$, $\ij$, is
continuous and of polynomial growth then using the result by
El-Karoui et al.(\cite{Elkarouietal97}, Thm. 8.5) related to connection between
solutions of reflected BSDEs and viscosity solutions of PDE with
obstacles we deduce that $(e^{-\l t}v^i(t,x))_{i\in \cJ}$ is a
viscosity solution for system (\ref{sysvi1}). \qed \bs

As previously mentioned, we now consider the case when the functions
$-f_i$, $\ij$, verify (H2)-(iv). Then we have:

\begin{theo}\label{exist_decreasingcase} Assume that assumptions (H3), (H4) are fulfilled and
that the functions $(-f_i)_{\ij}$ verify (H2). Then the system of
variational inequalities with inter-connected obstacles
(\ref{sysvi1}) has a continuous solution $(v^1,\dots,v^m)$ in the
class $\Pi^g$.
\end{theo}
\textbf{Proof:} We first prove existence of a candidate to be a
viscosity solution of the system. As previously, we will relate it
to the unique solution of the multidimensional reflected BSDE
(\ref{system}).\ms

\no {\underline{Step 1}: Construction}

For $\ij$, let $F_i$ be defined as in (\ref{expgrdf}). Choosing
$\lambda$ large enough, we obtain that each $F^{i}$ is decreasing
with respect to all the variables $y_{j}$, $j=1, \cdots, m$. Next
let us consider the following iterative Picard scheme: for any
$\ij$, $Y^{i, 0} = 0$ and for $n\geq 1$, we define $(Y^{i,n})_{\ij}$
by:
$$(Y^{1,n},\dots, Y^{m,n})=\Phi((Y^{1,n-1},\dots, Y^{m,n-1})).$$
Therefore we know from Theorem \ref{chassegneux} that the sequence
$(Y^{i,n})_{\ij}$ converges in $\cH^{2,m}$ to the unique solution
$(Y^{i})_{\ij}$ of the system of reflected BSDEs with oblique
reflection associated with
\\$((F_i(s,X^{t,x}_s,y^1,...,y^m,z))_{i\in \cJ}, (e^{\l
T}h_i(X^{t,x}_T))_{i\in \cJ}, (e^{\l t}g_{ij}(s,X^{t,x}_s))_{i,j\in
\cJ})$.

Next using Proposition \ref{sansz} and an induction argument, we deduce the
existence of continuous functions with polynomial growth such that:
$$\forall n\geq 0,
\forall \ij, \forall (t,x)\in \esp,\forall s\in
[t,T],\,\,Y^{i,n}_s=v^{i,n}(s,X^{t,x}_s).
$$
But from Theorem \ref{theo3}-(a) we have: for any $i,n,m$ and $s\leq
T$, \be\ba{l} \label{inegyi}\E[|Y^{i,n}_s-Y^{i,m}_s|^2]\leq C
\|(Y^{i,n-1})_{\ij}-(Y^{i,m-1})_{\ij}\|_{{\cH}^{2,m}}^{2} \mbox{ and }\\
\E[|Y^{i}_s-Y^{i,m}_s|^2]\leq C
\|(Y^{i})_{\ij}-(Y^{i,m-1})_{\ij}\|_{{\cH}^{2,m}}^{2}, \ea\ee where, to
obtain the last inequality, we rely on the characterization of the
solution $(Y^{i})_{\ij}$ constructed in \cite{chassagneuxetal2010},
i.e., $ (Y^{i})_{\ij}= \Phi((Y^{i})_{\ij})$. Taking $s=t$ we obtain:
$\forall \,\,(t,x)\in \esp$,
$$ \E[|Y_{t}^{i,n}-Y_t^{i,m}|] =  |v^{i,n}(t,x)-v^{i,m}(t,x)|\leq C \|(Y^{i,n-1})_{\ij}-(Y^{i,m-1})_{\ij}\|_{{\cH}^{2,m}}.
$$
As the sequence $((Y^{i,n})_{\ij})_{n\geq 0}$ is convergent in
${\cH}^{2,m}$ then it is of Cauchy type which implies that
$(v^{i,n})_{n\geq 0}$ is so and then converges pointwisely to a
deterministic function $v^i$, for any $\ij$. Thus going back to
(\ref{inegyi}) we deduce that: \be \label{representy}\forall \ij,\,\forall s\in [t,T],\, \P-a.s.,\,\,
Y^{i}_s=v^i(s,X^{t,x}_s).\ee
Let us now show that $v^i$, $\ij$, belongs to $\Pi^g$.
 Actually since $\Phi$ is a contraction in $({\cH}^{2,m},\|.\|_{\beta_0})$ and by some induction procedure on $n$ we get:
$$
\forall n,q\geq 0, \|(Y^{i,n+q})_{\ij}-(Y^{i,n})_{\ij}\|_{\beta
_0}\leq \frac{C_{\Phi}^n}{1-C_{\Phi}}\|(Y^{i,1})_{\ij}\|_{\beta
_0}\,;
$$
with $C_{\Phi}$ (such that $0< C_{\Phi} <1$) which is the contraction constant of the mapping $\Phi$ (constant which is independent of $(t,x)$).
As the norms $\|.\|$ and $\|.\|_{\beta_0}$ are equivalent then there
exists a constant $C_1$ such that:
$$
\forall n,q\geq 0, \|(Y^{i,n+q})_{\ij}-(Y^{i,n})_{\ij}\|\leq
C_1{C^n_{\Phi}}\|(Y^{i,1})_{\ij}\|.$$ Taking now the limit as $q$
goes to $+\infty$ and in view of (\ref{inegyi}) and
(\ref{representy}), if we then take $s=t$ we deduce that:
$$\forall (t,x)\in \esp,\,\,|v^i(t,x)-v^{i,n}(t,x)|\leq C_2 \|(Y^{i,1})_{\ij}\|.
$$
Finally one can check easily that $\|(Y^{i,1})_{\ij}\|(t,x)$ is of polynomial growth
(since $\E[\sup_{s\leq T}|X^{t,x}_s|^\gamma]$ belongs to $\Pi^g$ for any $\gamma \geq 0$, see (\ref{croix})) and since $v^{i,n}$ is so, then we deduce that $v^i$ is also of polynomial growth for
any $\ij$.
\medskip

\noindent \underline{Step 2}: Continuity of $v^i$, $\ij$.
\\
We again rely on the convergence result of any sequence $(Y^{i,n})_{n}$ constructed via the Picard iterative scheme. So let us initialize the scheme as follows:
 \begin{equation}\label{polynomgrowth}
\forall \; i \in \; \{1,\cdots,m\}, \forall \; s\le T, \quad   Y_{s}^{i, 0} = C\left( 1 + |X_{s}^{t,x}|^{p}\right),
 \end{equation}
 with the constant $C$ and the integer $p$ given by the fact that $ {v}^{i}$ is in $\Pi^{g}$.

 Next let $({Y}^{i, t, x})_{i =\{1,\cdots, m\}} $ be the unique solution of the
 multidimensional RBSDE associated with
$((F_i(s,X^{t,x}_s,y^1,...,y^m,z))_{i\in \cJ}, (e^{\l
T}h_i(X^{t,x}_T))_{i\in \cJ}, (e^{\l t}g_{ij}(s,X^{t,x}_s))_{i,j\in
\cJ})$ then we have
 $
 {Y}_{s}^{i, t, x} = {v}^{i}(s,\; X_{s}^{t,x})
 $ for any $s\in [t,T]$ and $\ij$. By definition of $Y^{i,0}$ in \ref{polynomgrowth} and the hypothesis on $v^{i}$, we obtain
 $$
 \P-a.s.,\forall \ij,\forall s\in [t,T],  {Y}_{s}^{i, t, x}\leq Y_{s}^{i, 0}.$$
 Let us now prove by induction that for any $n\geq 0$ we have:
\begin{equation}\label{eq:lower_and_upperestimates}
\dis{ \forall \; n \in \mathbb{N}, \quad \forall \ij,\, \forall s\in [t,T],\,\, Y_{s}^{i, 2n+1} \le {Y}_{s}^{i}
 \le Y_{s}^{i, 2n}.}\end{equation}
For $n=0$, the inequality of the right-hand side is already true. Let consider the left-hand side.
Recalling that by construction we have $ (Y^{i,1})_{\ij} = \Phi((Y^{i,0})_{\ij})$ and $({Y}^{i})_{\ij} = \Phi(({Y}^{i})_{\ij})$,
we get the following comparison result for the two drivers of the RBSDE satisfied by ${Y}^{i}$ and $Y^{i,1}$
 $$ F_{i}(s, X^{t,x}_{s}, {Y}_{s}^{1},\cdots, {Y}_{s}^{m} ,z^{i} ) \ge F_{i}(s, X^{t,x}_{s}, Y_{s}^{1,0}, \cdots , Y_{s}^{m,0} , z^{i} ),$$since $F_i$ is non-increasing w.r.t. all variables $y_{j}$, $j\in \cJ$. Next using the comparison result
 (see Remark \ref{comparison}) and uniqueness
of the solution of the system associated with \\ $((F_i(s,X^{t,x}_s,y^1,...,y^m,z))_{i\in \cJ}, (e^{\l
T}h_i(X^{t,x}_T))_{i\in \cJ}, (e^{\l t}g_{ij}(s,X^{t,x}_s))_{i,j\in
\cJ})$ we deduce that $Y^{i,1}\leq Y^i$ for any $\ij$. Thus the property (\ref{eq:lower_and_upperestimates}) is valid for $n=0$. Now if it is satisfied for some $n$ and repeating the same argumentation, it also holds for $n+1$, whence the claim.

Next and relying once more on the result obtained 
El-Karoui et al.(\cite{Elkarouietal97}, Thm. 8.5) let $\tilde v^{i,n}$, $\ij$ and $n\geq 0$, be the deterministic continuous functions of $\Pi^g$ such that:$$
\forall (t,x)\in \esp, \forall s\in [t,T], Y^{i,n}_s=\tilde v^{i,n}(s,X^{t,x}_s).
$$
First and as previously for any $i\in \cJ$, the sequence
$(\tilde v^{i,n})_{n\geq 0}$ is of Cauchy type and converges
pointwisely to $v^i$. Next the inequalities
(\ref{eq:lower_and_upperestimates}) imply that: $\forall n\geq 0$
and $\ij$,
$$\tilde v^{i,2n+1}\leq v^i\leq \tilde v^{i,2n},
$$
which implies 
$$  v^i = \lim_{n} \nearrow \tilde v^{i,2n+1} = \lim_{n} \searrow \tilde v^{i,2n}. $$
Therefore for any $\ij$, $v^i$ is both $usc$ and $lsc$ and thus it is continuous. Next, since $(Y^i,Z^i,K^i)_{\ij}$ is the unique
solution of the system of reflected BSDEs associated with the following triplet of datas \\
$((F_i(s,X^{t,x}_s,y^1,...,y^m,z))_{i\in \cJ}, (e^{\l
T}h_i(X^{t,x}_T))_{i\in \cJ}, (e^{\l t}g_{ij}(s,X^{t,x}_s))_{i,j\in
\cJ})$ then\\ $((e^{-\l t}Y^i_t,e^{-\l t}Z^i_t,e^{-\l
t}dK^i_t)_{t\leq T})_{\ij}$ is the solution of the sytem of
reflected BSDEs associated with
$((f_i)_{\ij},(h_i)_{\ij},(g_{ij})_{i,j\in \cJ})$ then using once
more the result by El-Karoui et al. (\cite{elkarouiquenezpeng},
Thm.8.5) to deduce that $(e^{-\l t}v^i)_{\ij}$ is a continuous with
polynomial growth solution of the system of variational inequalities
with inter-connected obstacles (\ref{sysvi1}). The proof is now
complete. \qed \ms

Next we deal with the issue of uniqueness of the solution of
(\ref{sysvi1}) in the general case.
\begin{theo}
Under Assumptions (H1),(H2),(H3) and (H4), the solution of the
system of variational inequalities with inter-connected obstacles
(\ref{sysvi1}) is unique in the class $(v_1,....,v_m)$ of continuous
functions which belong to $\Pi^g$.
\end{theo}
\underline{Proof}. As usual it is enough to show that if
$(u_1,...,u_m)$ (resp. $(v_1,...,v_m)$) is a continuous subsolution
(resp. supersolution) of (\ref{sysvi1}) such that $u_i,v^i$, $\ij$,
belong to $\Pi^g$ then $u_i\le v^i$, for all $\ij$. Thus classically
we have uniqueness of (\ref{sysvi1}). \ms

\no \underline{Step 1}: We first assume the existence of a constant
$\l<-m.\max \{C^j_f, j=1,...,m\}$ ($C^j_f$ is the Lipschitz constant
of $f_j$ involved in [H2]-(ii)) such that for any $\ij$, $f_i$
verifies:

$\forall$ $t,x,y_1,...,y_{i-1},y_{i+1},...,y_m, y, \bar y,$ if
$y\geq \bar y$ then \be\label{hypinter}\begin{array}{l}
f_i(t,x,y_1,...,y_{i-1},y,y_{i+1},...,y_m)-f_i(t,x,y_1,...,y_{i-1},\bar
y,y_{i+1},...,y_m)\leq \l (y-\bar y).\end{array}\ee

Let $\gamma>0$ and $C$ be such that for any $i\in \cJ$ we have:
$$
|u_i(t,x)|+|v^i(t,x)|\leq C(1+|x|^\g),\,\,\forall (t,x)\in \esp.$$

Next as in Lemma \ref{modifsursolution}, for $\theta
>0$ and $\n$ large enough $(v^i(t,x)+\theta e^{-\n t}|x|^{2\g +2})_{i=1,m}$ is also a
supersolution for (\ref{system}). Therefore it is enough to show
that for any $i\in \cJ$, we have:
$$
\forall (t,x)\in \esp, u_i(t,x)\leq v^i(t,x)+\theta e^{-\n
t}|x|^{2\g +2},$$ and taking the limit as $\theta \rightarrow 0$ we
obtain the desired result. So let us set
$w^{i,\theta,\n}(t,x)=v^i(t,x)+\theta e^{-\n t}|x|^{2\g +2}$,
$(t,x)\in \esp$ which is denoted by  $w^i$ for simplicity. Assume
now that there exists a point $(\bar t,\bar x)\in \esp$ such that
$\max_{i\in \cJ}(u^i(\bar t,\bar x)-w^i(\bar t,\bar x))>0.$ Using
the growth condition on $u_{i}$ and $w^{i}$, there exists $R>0$ such
that:
$$\forall \ij,\,\,
\forall (t,x)\in \esp \mbox{ s.t. }|x|\geq R,\,\,u_i(
t,x)-w^i(t,x)<0.$$ Taking into account the values of the subsolution
and the supersolution at $T$, it implies that
\be\label{sgn1}\begin{array}{l} 0<\max_{(t,x)\in \esp }\max_{i\in
\cJ}(u_i( t,x)-w^i(t,x))=\\\qquad \qq \max_{(t,x)\in [0,T[\times
B(0,R)}\max_{i\in \cJ}(u_i( t,x)-w^i(t,x))=\max_{i\in \cJ}(u_i(
t^*,x^*)-w^i(t^*,x^*))\end{array} \ee where $(t^*,x^*)\in
[0,T[\times B(0,R)$. Now let $\tilde \cJ$ be as in (\ref{tildej})
and let $j\in \tilde \cJ$ be such that \be\label{penet}
u_j(t^*,x^*)>\max_{k\in \cJ^{-j}}(u_k(t^*,x^*)-g_{jk}(t^*,x^*)).\ee
Next let $\Phi^{j}_n(t,x,y)$ and $\varphi_n$, $n\geq 0$, be the same
functions as in (\ref{eqphin}) and let $(t_n,x_n,y_n)\in [0,T]\times
B'(0,R)^2$ be the triple such that:
$$
\Phi^{j}_n(t_n,x_n,y_n)=\max_{(t,x,y)\in [0,T]\times
B'(0,R)^2}\Phi^{j}_n(t,x,y).
$$
As in the proof of Theorem \ref{comparisonviscosity}, one can show
that: \be \label{convxn}(t_{n},x_{n},y_{n}) \rw (t^*,x^*,x^*) \mbox{
and }n|x_n-y_n|^{2\g+2}\rightarrow 0 \mbox{ as }n\rw \infty.\ee Now
as $u_j$ and $g_{jk}$ are continuous functions and taking into
account (\ref{penet}) we deduce that for $n$ large enough:
\be\label{cil2x}u_j(t_n,x_n)>\max_{k\in
\cJ^{-j}}(u_k(t_n,x_n)-g_{jk}(t_n,x_n)).\ee We next apply
Crandall-Ishii-Lions's Lemma (see e.g. \cite{Crandalllions92} or
\cite{Flemingetsoner}, pp.216) with $\Phi_n^{j}$, $u_j$, $w^j$ and
$\varphi_n$ (recall that (\ref{cil2x}) is satisfied) at
$(t_n,x_n,y_n)$, there exist $(p^n_u,q^n_u,M^n_u)\in \bar
J^{2,+}(u_j)(t_n,x_n)$ and $( p^n_w,q^n_w,M^n_w)\in \bar
J^{2,-}(w^j)(t_n,y_n)$ such that:

$p^n_u-p^n_w=\partial_t\tilde \varphi_n(t_n,x_n,y_n)=2(t_n-t^*)$,
$q^n_u \,\,(\mbox{resp. }q^n_w)\,=\partial_x
\varphi_n(t_n,x_n,y_n)$$(\mbox{resp.}-\partial_y
\varphi_n(t_n,x_n,y_n)) \mbox{ and }$ \be \label{cil31}
\left (\begin{array}{ll} M_u^n&0\\
0&-N_w^n\end{array}\right )\leq A_n+\frac{1}{2n}A_n^2\ee where
$A_n=D^2_{(x,y)}\varphi_n(t_n,x_n,y_n)$. Taking now into account
that $(u_i)_{i\in \cJ}$ (resp. $(w^i)_{i\in \cJ}$) is a subsolution
(resp. supersolution) of (\ref{sysvi1}) and once more (\ref{cil2x})
we deduce:
$$
-p^n_u-b(t_n,x_n)^\top.q^n_u-\frac{1}{2}
Tr[(\sigma\sigma^\top)(t_n,x_n)M^n_u]-f_j(t_n,x_n,(u_i(t_n,x_n))_{\ij},\sigma
(t_n,x_n)^\top.q^n_u)\leq 0
$$
and
$$-p^n_w-b(t_n,y_n)^\top.q^n_w-\frac{1}{2}
Tr[(\sigma\sigma^\top)(t_n,y_n)M^n_w]-f_j(t_n,y_n,(w^i(t_n,y_n))_{\ij},\sigma
(t_n,y_n)^\top.q^n_w)\geq 0.
$$
Making the difference between those two inequalities yields:
\be\label{cil42}\begin{array}{l}
-(p^n_u-p^n_w)-(b(t_n,x_n)^\top.q^n_u-b(t_n,y_n)^\top.q^n_w)-
\frac{1}{2}
Tr[\{\sigma\sigma^\top(t_n,x_n)M^n_u-\sigma\sigma^\top(t_n,y_n)M^n_w\}]\\\qquad
-\{f_j(t_n,x_n,(u_i(t_n,x_n))_{\ij},\sigma (t_n,x_n)^\top.q^n_u)
-f_j(t_n,y_n,(w^i(t_n,y_n))_{\ij},\sigma (t_n,y_n)^\top.q^n_w)\}\leq
0 \end{array} \ee and then
$$\begin{array}{l}
-\{f_j(t_n,x_n,(u_i(t_n,x_n))_{\ij},\sigma (t_n,x_n)^\top.q^n_u)
-f_j(t_n,x_n,(w^i(t_n,y_n))_{\ij},\sigma (t_n,x_n)^\top.q^n_u))\leq
A_n\end{array}
$$
where as in (\ref{estibn}) and (\ref{estisn}) we have $\limsup_{n\rw
\infty}A_n\le 0$. Next linearizing $f_j$, which is Lipschitz w.r.t.
$(y^i)_{\ij}$, and using (\ref{hypinter}) we obtain:
$$
-\l
(u_j(t_n,x_n)-w^j(t_n,y_n))-\sum_{k\in\cJ^{-j}}\Theta_n^{j,k}((u_k(t_n,x_n)-w^k(t_n,y_n))\leq
A_n
$$
with $\Theta_n^{j,k}$ which stands for the increment rate of $f_j$ w.r.t. $y_k$ and
which, thanks to the monotonicity assumption (see [H2] (iv)), is nonnegative and bounded by $C_f^j$ the Lipschitz constant
of $f_j$. Thus
$$\ba{ll}
-\l (u_j(t_n,x_n)-w^j(t_n,y_n))&\leq
\sum_{k\in\cJ^{-j}}\Theta_n^{j,k}((u_k(t_n,x_n)-w^k(t_n,y_n))^++
A_n\\
{}&\leq C_f^j\sum_{k\in\cJ^{-j}}((u_k(t_n,x_n)-w^k(t_n,y_n))^++ A_n.
\ea
$$
Taking now the limit as $n\rw \infty$ and since $j\in \tilde {\cal
J}$ to obtain:
$$\ba{ll}
-\l(u_j(t^*,x^*)-w^j(t^*,x^*))&\leq C_f^j\sum_{k\in
\cJ^{-j}}(u_k(t^*,x^*)-w^k(t^*,x^*))^+\\{}&\leq (m-1) C_f^j
(u_j(t^*,x^*)-w^j(t^*,x^*))\ea
$$
which contradictory since $u_j(t^*,x^*)-w^j(t^*,x^*)>0$ and $-\l> m
C^j_f$. Thus for any $\ij$, $u_i\leq w^i$. \ms

\noindent \underline{Step 2}: The general case.

For arbitrary $\l\in \R$, if $(u_j)_{j\in \cJ}$ (resp. $(v^j)_{j\in
\cJ}$) is a subsolution (resp. supersolution) of (\ref{sysvi1}) then
$\tilde u_j(t,x)=e^{\l t}u_j(t,x)$ and $\tilde v^j(t,x)=e^{\l
t}v^j(t,x)$ is a subsolution (resp. supersolution) of the following
system of variational inequalities with oblique reflection (see the
proof of Proposition \ref{comparisonviscosity}, Step 2): $\forall
\,\,i\in \cJ$,\be \label{sysvi3} \left\{
\begin{array}{l}
\min\left \{\tilde v_i(t,x)- \max\limits_{j\in{\cal J}^{-i}}(-e^{\l
t}g_{ij}(t,x)+\tilde v_j(t,x))\right.,\\ \qquad \qquad
\left.-\partial_t\tilde v_i(t,x)- {\cal
L}\tilde v_i(t,x)+\l \tilde v_i(t,x)-e^{\l t}f_i(t,x,(e^{-\l t}
\tilde v^i(t,x))_{\ij},e^{-\l t}\sigma^\top(t,x).D_x \tilde v^i(t,x))\right\}=0;\\
\tilde v_i(T,x)=e^{\l T}h_i(x).
\end{array}\right.
\ee But in choosing $\l$ small enough the functions
$F_i(t,x,(u_i)_{\ij},z)=-\l u+e^{\l t}f_i(t,x,(e^{-\l
t}u_i)_{\ij},e^{-\l t}z)$, $\ij$, satisfy condition (\ref{hypinter})
and then thanks to the result stated in Step 1, we have $\tilde
u_i\le \tilde v_i$, $\ij$. Thus $u_i\le v_i$ for any $\ij$ which is
the desired result.  \qed \ms

As a by-product we have:

\begin{corollary} Assume that (H1), (H3), (H4) and (H2) are
fulfilled, then there is a unique solution of the system
(\ref{sysvi1}) in the class of continuous functions with polynomial
growth. \qed
\end{corollary}

\small{
}
\end{document}